\documentclass{article}
\usepackage{amsmath}
\usepackage{multicol,color}
\usepackage{float}
\usepackage{soul}
\usepackage{graphicx}
\usepackage{amssymb}
\usepackage{theorem}
\usepackage{fancyhdr}
\usepackage{tikz}
\usepackage{enumerate}
\usepackage[margin=1.1in]{geometry}
\usepackage{pgfplots}
\usepgfplotslibrary{polar}
\usepgflibrary{shapes.geometric}
\pgfplotsset{my style/.append style={axis x line=middle, axis y line= middle, xlabel={$x$}, ylabel={$y$}, axis equal }}
\usetikzlibrary{through,calc}
\usetikzlibrary{calc,intersections} 
\newcommand\markangle[6][red]{
  \begin{scope}
    \path[clip] (#2) -- (#3) -- (#4);
    \fill[color=#1,fill opacity=0.5,draw=#1,name path=circle]
    (#3) circle (#6mm);
  \end{scope}
  \path[name path=line one] (#3) -- (#2);
  \path[name path=line two] (#3) -- (#4);
  \path[%
  name intersections={of=line one and circle, by={inter one}},
  name intersections={of=line two and circle, by={inter two}}
  ] (inter one) -- (inter two) coordinate[pos=.5] (middle);
  \path[%
  name path=bissectrice
  ] (#3) -- (barycentric cs:#3=-1,middle=1.2);
  \path[
  name intersections={of=bissectrice and circle, by={middleArc}}
  ] (#3) -- (middleArc) node[pos=1.3] {#5};
  }

\def\disp{\displaystyle}
\def\ve{\varepsilon}

\def\lm{\lambda}
\def\O{\Omega}

\def\oa{\bar a}

\def\ox{\bar{x}}

\def\oz{\bar{z}}

\def\ou{\bar{u}}

\def\gph{\hbox{}}

\def\gg{\gamma}

\def\Lto{\Longrightarrow}
\def\tto{\rightrightarrows}

\def\Bar{\overline}
\def\ra{\rangle}
\def\la{\langle}
\def\ve{\varepsilon}

\def\h{\hfill\Box}
\def\R{\mathbb{R}}
\def\N{\mathbb{N}}

\def\co{\mbox{\rm co}\,}
\def\gph{\mbox{\rm gph}\,}
\def\epi{\mbox{\rm epi}\,}

\def\dist{\mbox{\rm dist}}

\def\bd{\mbox{\rm bd}\,}

\def\O{\Omega}

\def\ph{\varphi}
\def\emp{\emptyset}

\def\lm{\lambda}

\def\gg{\gamma}

\def\al{\alpha}

\def\be{\beta}

\def\N{I\!\!N}
\def\th{\theta}
\def\vTh{\vartheta}
\newtheorem{theorem}{Theorem}[section]

\newtheorem{corollary}[theorem]{Corollary}
\newtheorem{proposition}[theorem]{Proposition}
\newtheorem{definition}[theorem]{Definition}
\newtheorem{example}[theorem]{Example}
\newtheorem{remark}[theorem]{Remark}

\begin{document}
\begin{center}
\vspace*{0.3in}{\bf APPLICATIONS OF OPTIMAL CONTROL OF A NONCONVEX SWEEPING PROCESS TO OPTIMIZATION OF THE PLANAR CROWD MOTION MODEL}\\[3ex]
TAN H. CAO\footnote{Department of Applied Mathematics and Statistics, State University of New York--Korea, Yeonsu-Gu, Incheon, Republic of Korea (tan.cao@stonybrook.edu). Research of this author was partly supported by the Ministry of Science, ICT and Future Planning, Republic of Korea under the ICT Consilience Creative Program (IITP-2015-R0346-15-1007) supervised by the Institute for Information \& Communications Technology Promotion.}
and B. S. MORDUKHOVICH\footnote{Department of Mathematics, Wayne State University, Detroit, Michigan, USA and RUDN University, Moscow 117198, Russia (boris@math.wayne.edu). Research of this author was partly supported by the US National Science Foundation under grant DMS-1512846, by the US Air Force Office of Scientific Research under grant \#15RT0462, and by the RUDN University Program 5-100.}\\[2ex]
{\bf Dedicated to Peter E. Kloeden on occasion of his 70th birthday}\\[2ex]
\end{center}
\small{\bf Abstract.} This paper concerns optimal control of a nonconvex perturbed sweeping process and its applications to optimization of the planar crowd motion model of traffic equilibria. The obtained theoretical results allow us to investigate a dynamic optimization problem for the microscopic planar crown motion model with finitely many participants and completely solve it analytically in the case of two participants.\\[1ex]
\vspace*{-0.1in}
\normalsize

\section{Introduction}\setcounter{equation}{0}
\vspace*{-0.1in}

This paper is mainly devoted to applications of very recent and new results on optimal control of a perturbed nonconvex sweeping process to solving dynamic optimization problems for a microscopic version of the planar crowd motion model the importance of which has been well recognized in traffic equilibria and other areas of socioeconomic modeling; see, e.g., \cite{mv2,ve} and the references therein. Our previous attempts in this direction \cite{cm1,cm2} were limited by the crowd motion model in a {\em corridor} whose optimal control description is a particular case of the following dynamic optimization problem governed by the convex {\em controlled polyhedral sweeping process}:
\begin{eqnarray}\label{e:Bolza}
\mbox{minimize }\;J[x,u,a]:=\ph\big(x(T)\big)+\int^T_0\ell\big(t,x(t),u(t),a(t),\dot{x}(t),\dot{u}(t),\dot{a}(t)\big)dt
\end{eqnarray}
over the control functions $u(\cdot)\in W^{1,2}([0,T];\R^n)$ and $a(\cdot)\in W^{1,2}([0,T];\R^d)$ and the corresponding trajectories $x(\cdot)\in W^{1,2}([0,T];\R^n)$ of the differential inclusion
\begin{eqnarray}\label{e:psp}
-\dot{x}(t)\in N\big(x(t);C(t)\big)+f\big(x(t),a(t)\big) \ \mbox{ a.e. } \ t\in [0,T], \ x(0):=x_0\in C(0),
\end{eqnarray}
where $x_0\in\R^n$ and $T>0$ are fixed and where the controlled moving set is given by
\begin{eqnarray}\label{poly}
C(t):=C+u(t)\;\mbox{ with }\;C:=\big\{x\in\R^n|\;\la x^*_i,x\ra\le 0\;\mbox{ for all }\;i=1,\ldots,m\big\}.
\end{eqnarray}
In the description above, the differential inclusion $-\dot{x}(t)\in N(x(t);C(t))$ signifies the basic Moreau's {\em sweeping process} \cite{mor_frict} introduced and studied in the case of a {\em uncontrolled} nicely moving convex set $C(t)$ with the symbol $N(x;C)$ standing the classical {\em normal cone} to a convex set $C$ at $x$. It has been well realized in the sweeping process theory that the Cauchy problem for the latter inclusion with $x(0)=x_0$ has a {\em unique} solution, and hence does not give us any room for optimization.\vspace*{-0.05in}

It seems that for the first time in the literature optimal control problems for the sweeping process with {\em controlled} sweeping sets $C(t)$ as in \eqref{poly} were formulated and studied in the papers by Colombo et al.\ \cite{chhm1,chhm2}, where necessary optimality conditions for such problems were derived. However, the aforementioned papers did not address the sweeping dynamics \eqref{e:psp} with the other type of controlled functions $a(\cdot)$ entering the external force $f$ in \eqref{e:psp}. Such problems \eqref{e:Bolza}--\eqref{poly} were first considered in \cite{cm1,cm2} with applying the necessary optimality conditions derived therein to solving optimal control problems for the corridor version of the crowd motion model formulated in this form. The corridor restriction of the crowd motion model investigated in \cite{cm1,cm2} was due to the fact that only this version fitted the polyhedral description of $C$ in \eqref{poly}, which was essentially exploited in \cite{cm1,cm2} as well as in \cite{chhm1,chhm2}. The natural desire to consider the much more practical {\em planar} version of the crowd motion model mandatory requires to deal with {\em nonpolyhedral} (actually {\em nonconvex}) sweeping sets as shown in \cite{mv2,ve}.\vspace*{-0.05in}

Before proceeding in this direction, note that even optimal control problems of type \eqref{e:Bolza}--\eqref{poly} are significantly different from and essentially more involved than those in standard control theory and dynamic optimization governed by {\em Lipschitzian} differential inclusions considered, e.g., in \cite{clsw,m-book2,v}. Discontinuity of the normal cone mapping $x\mapsto N(x;C)$ is the main obstacle to employ various approaches developed in optimal control of differential inclusions. More complications come from the intrinsic presence of pointwise state constraints of equality and inequality types in \eqref{e:psp} since we have $N(x;C)=\emp$ for $x\notin C$. These challenges were met in \cite{cm1,cm2,chhm1,chhm2} for the polyhedral processes \eqref{e:Bolza}--\eqref{poly} by extending the {\em method of discrete approximations} developed in \cite{m95,m-book1} for optimal control of Lipschitzian differential inclusions.\vspace*{-0.05in}

The origin of using the time discretization and finite-difference approximations to derive optimality conditions for infinite-dimensional variational problems goes back to Euler at the very beginning of the calculus of variations. However, nowadays Euler's and related finite-difference schemes are phenomenally used in numerical methods of solving various problems for deterministic and stochastic dynamic systems; see, e.g., the now classical book \cite{klo-pat92} and other publications by Peter Kloeden and his collaborators. In contrast, our approach---while having some numerical flavor---employs Euler's finite differences mainly as a vehicle to obtain necessary optimality conditions for the nonstandard classes of optimal control problems under consideration and their valuable applications. In this way we combine the method of discrete approximations with deep developments in variational analysis and generalized differentiation.\vspace*{-0.05in}

Note that some optimal control problems for several versions of the sweeping process, different from \eqref{e:Bolza}--\eqref{poly} and the one studied below, have been recently considered in \cite{ao,ac,bk,pfs} with deriving necessary optimality conditions by using other approximation techniques. However, the results obtained therein cannot be applied to either corridor or planar forms of the crowd motion model.\vspace*{-0.05in}

An appropriate framework of sweeping optimal control to encompass optimization of the planar crowd motion model of our interest in this paper is described as minimizing the generalized Bolza functional \eqref{e:Bolza} over the controlled perturbed sweeping process \eqref{e:psp} with $u(\cdot)\in W^{1,2}([0,T];\R^n)$, $a(\cdot)\in W^{1,2}([0,T];\R^d)$, and $x(\cdot)\in W^{1,2}([0,T];\R^n)$ as above, but with replacing the convex polyhedral sweeping set \eqref{poly} by the {\em nonconvex} (and hence nonpolyhedral) one
\begin{eqnarray}\label{e:mset}
C(t):=C+u(t)=\bigcap^m_{i=1}C_i+u(t)\;\mbox{ with }\;C_i:=\big\{x\in\R^n\big|\;g_i(x)\ge 0\big\},\quad i=1,\ldots,m,
\end{eqnarray}
defined by some convex and ${\cal C}^2$-smooth functions $g_i\colon\R^n\to\R$. Since the set $C(t)$ is nonconvex for any $t\in[0,T]$, we have to specify the suitable notion of the normal cone in \eqref{e:psp}. It occurs, however, that all the major normal cone of variational analysis agree with each other under the assumptions imposed on $g_i$; see Section~2. Thus we can keep the generic normal cone symbol ``$N$" in our setting.\vspace*{-0.05in}

Besides the dynamic constraints in \eqref{e:psp} and \eqref{e:mset}, we impose the following pointwise constraints on the control functions $u(t)$ in \eqref{e:mset}:
\begin{eqnarray}\label{e:ec}
r_1\le\|u(t)\|\le r_2\;\mbox{ for all }\;t\in[0,T]
\end{eqnarray}
with the fixed bounds $0<r_1\le\|u(t)\|\le r_2$, where it might be that $r_1=r_2$, and thus \eqref{e:ec} becomes a pointwise constraint of the equality type. The other type of pointwise constraints on control and state functions intrinsically arise from \eqref{e:psp}, \eqref{e:mset} due to $N(x(t);C(t))=\emp$ for $x(t)\notin C(t)$ and are written as
\begin{eqnarray}\label{e:mc}
g_i\big(x(t)-u(t)\big)\ge 0\;\mbox{ for all }\;t\in[0,T],\quad i=1,\ldots,m.
\end{eqnarray}\vspace*{-0.2in}

Throughout the paper, {\em denote by $(P)$} the constrained optimal control problem for a nonconvex sweeping process described in \eqref{e:Bolza}, \eqref{e:psp}, \eqref{e:mset}, \eqref{e:ec}, and \eqref{e:mc}. We can equivalently reformulate problem $(P)$ in the conventional form of dynamic optimization for differential inclusions with respect to the triple $z:=(x,u,a)\in\R^n\times\R^n\times\R^d$ over $z(\cdot)\in W^{1,2}([0,T];\R^{2n+d})$ satisfying the conditions
\begin{eqnarray}\label{e:psp1}
-\dot z(t)\in F\big(z(t)\big)\times\R^n\times\R^d\;\mbox{ a.e. }\;t\in[0,T],\quad x_0-u(0)\in C\;\mbox{ for }\;z(0)=\big(x_0,u(0),a(0)\big),
\end{eqnarray}
with $F(z):=N(x-u;C)$ and $C$ taken from \eqref{e:psp}, under the pointwise {\em state constraints} in \eqref{e:ec} and \eqref{e:mc} imposed on trajectories $z(t)$ of \eqref{e:psp1} on the whole interval $[0,T]$. Note however that the mapping
\begin{eqnarray*}
G(z)=G(x,u,a):=N(x-u;C)\times\R^n\times\R^d,\quad z\in\R^{2n+d},
\end{eqnarray*}
is highly nonstandard in control theory for differential inclusions while being intrinsically {\em non-Lipschitzian} (in fact discontinuous) and unbounded. Nevertheless, an appropriate development of the method of discrete approximations combined with advanced tools of first-order and second-order variational analysis and generalized differentiation allows us to derive necessary optimality conditions for $(P)$. As a by-product of this approach, we establish the strong $W^{1,2}$-convergence of optimal solutions for finite-dimensional discrete-time approximation problems to the reference optimal solution for the sweeping control problem $(P)$. Thus finite-dimensional optimal solutions for discrete approximations can be treated as constructive {\em suboptimal} solutions for the original continuous-time problem under consideration.\vspace*{-0.05in}

This program has been basically realized in our recent preprint \cite{cm3}. For completeness and reader's convenience we present in Section~2 the main results of \cite{cm3}, which are applied here to the planar crowd motion model; it was our major motivation for \cite{cm3}. Then Section~3 contains applications of the obtained results for the sweeping optimal control to an optimization problem for the planar crowd motion model while Section~4 is devoted to solving the latter problem with two participants. In the final Section~5 we summarize the main developments and discuss some directions of out future research.\vspace*{-0.05in}

The notation of this paper is standard in variational analysis and optimal control; see, e.g., \cite{m-book1,rw,v}. Recall that the symbol $B(x,\ve)$ denotes the closed ball of the space in question centered at $x$ with radius $\ve>0$ while $\N$ signifies the collections of all natural numbers $\{1,2,\ldots\}$.\vspace*{-0.2in}

\section{Optimal Control of Nonpolyhedral Sweeping Processes}
\setcounter{equation}{0}\vspace*{-0.1in}

The purpose of this section is to overview major results of our preprint \cite{cm3} (available in arXiv) and adapt them for the subsequent applications to the planar crowd motion model. We first recall some notions of variational analysis and generalized differentiation that are systematically used in below.\vspace*{-0.05in}

Given a set $\O\subset\R^n$ locally closed around $\ox\in\O$, consider the distance function
\begin{eqnarray*}
\dist(x;\O):=\min\big\{\|x-y\|\;\big|\;y\in\O\big\}
\end{eqnarray*}
for $x\in\R^n$ near $\ox$ and define the Euclidean projector of $x$ onto $\O$ by
\begin{eqnarray*}
\Pi(x;\O):=\big\{w\in\O\big|\;\|x-w\|=\dist(x;\O)\big\}.
\end{eqnarray*}
Then the {\em proximal normal} cone to $\O$ at $\ox$ is given by
\begin{eqnarray}\label{pnc}
N^P(\ox;\O):=\big\{v\in\R^n\big|\;\exists\,\al>0\;\mbox{ such that }\;\ox\in\Pi(\ox+\al v;\O)\big\},\quad\ox\in\O,
\end{eqnarray}
with $N^P(\ox;\O):=\emp$ if $\ox\notin\O$. Using this construction allows us to formulate the notion of uniform prox-regularity for sets that plays a crucial role in the subsequent developments.\vspace*{-0.05in}

\begin{definition}{\bf (uniform prox-regularity of sets).}\label{Def-1} Let $\O$ be a closed subset of $\R^n$, and let $\eta>0$. Then $\O$ is $\eta$-{\sc uniformly prox-regular} if for all $x\in\bd\O$ and $v\in N^P(x;\O)$ with $\|v\|=1$ we have $B(x+\eta v,\eta)\cap\O=\emp$. Equivalently, the $\eta$-uniform prox-regularity of $\O$ can be defined via the estimate
\begin{eqnarray*}
\la v,y-x\ra\le\dfrac{\|v\|}{2\eta}\|y-x\|^2\;\mbox{ for all }\;y\in\O,\;x\in\bd\O,\;\mbox{ and }\;v\in N^P(x;\O).
\end{eqnarray*}
\end{definition}\vspace*{-0.05in}

It is worth mentioning that convex sets are $\infty$-uniformly prox-regular and that for every $\eta>0$ the projection operator $\Pi(x;\O)$  for an $\eta$-uniformly prox-regular set is single-valued if $\dist(x;\O)<\eta$.\vspace*{-0.05in}

We refer the reader to the excellent survey \cite{CT} for many results and history of prox-regular sets (known in geometric measure theory as ``sets with positive reach") and its uniform version from Definition~\ref{Def-1}. It is important to observe that all the major normal cones of variational analysis (regular/Fr\'echet, limiting/Mordukhovich, and convexified/Clarke) agree with each other and with the proximal one for the uniformly prox-regular sets, and thus we use the symbol ``N" for all of them in \eqref{e:psp} and in what follows. In particular, we can freely use in our subsequent study and applications full calculi and computations available for the {\em limiting normal cone} to $\O$ at $\ox\in\O$ defined by
\begin{eqnarray}\label{e:Mor-nc}
N(\ox;\O):=\big\{v\in\R^n\big|\;\exists\,x_k\to\ox,\;w_k\in\Pi(x;\O),\;\al_k\ge 0\;\mbox{ s.t. }\;\al_k(x_k-w_k)\to v\;\mbox{ as }\;k\to\infty\big\}
\end{eqnarray}
as well as for the associated with it coderivative of set-valued mappings and the first-order and second-order subdifferentials of extended-real-valued functions.\vspace*{-0.05in}

Our {\em standing assumptions} on the initial data of $(P)$, which particularly ensure the uniform prox-regularity of the set $C(t)$ in \eqref{e:mset} for each $t\in[0,T]$, are the following:

{\bf (H1)} The perturbation mapping $f\colon\R^n\times\R^d\to\R^n$ in \eqref{e:psp} is continuous on $\R^n\times\R^d$ and locally Lipschitzian with respect to the first argument, i.e., for every $\ve>0$ there is a constant $K>0$ such that
\begin{eqnarray*}
\|f(x,a)-f(y,a)\|\le K\|x-y\|\;\mbox{ whenever }\;(x,y)\in B(0,\ve)\times B(0,\ve),\;a\in\R^d.
\end{eqnarray*}
Furthermore, there is a constant $M>0$ ensuring the growth condition
\begin{eqnarray*}
\|f(x,a)\|\le M\big(1+\|x\|\big)\;\mbox{ for any }\;x\in\bigcup_{t\in[0,T]}C(t),\;a\in\R^d.
\end{eqnarray*}\vspace*{-0.1in}

{\bf (H2)} There exist $c>0$ and open sets $V_i\supset C_i$ such that
\begin{eqnarray}\label{e:dist-inequality}
d_H(C_i;\R^n\setminus V_i)>c,
\end{eqnarray}
where $d_H$ denotes the Hausdorff distance between sets, and where $C_i$, $i=1,\ldots,m$, are taken in \eqref{e:mset}. In addition, there are positive constants $M_1,M_2$, and $M_3$ such that the functions $g_i(\cdot)$, $i=1,\ldots,m$, are ${\cal C}^2$-smooth satisfying the estimates
\begin{eqnarray}\label{e:bound-grad}
M_1\le\|\nabla g_i(x)\|\le M_2\;\mbox{ and }\;\|\nabla^2 g_i(x)\|\le M_3\;\mbox{ for all }\;x\in V_i.
\end{eqnarray}\vspace*{-0.1in}

{\bf (H3)} There exist positive numbers $\be$ and $\rho$ such that
\begin{eqnarray}\label{e:w-inverse-triangle-in}
\sum_{i\in I_\rho(x)}\lm_i\|\nabla g_i(x)\|\le\be\Big\|\sum_{i\in I_\rho(x)}\lm_i\nabla g_i(x)\Big\|\;\mbox{ for all }\;x\in C\;\mbox{ and }\;\lm_i\ge 0
\end{eqnarray}
with the index set $I_\rho(x):=\{i\in\{1,\ldots,m\}|\;g_i(x)\le\rho\}$.\vspace*{-0.05in}

{\bf (H4)} The terminal cost $\ph\colon\R^n\to\Bar{\R}$ is lower semicontinuous (l.s.c.), while the running cost $\ell$ in \eqref{e:Bolza} is such that $\ell_t:=\ell(t,\cdot)\colon\R^{4n+2d}\to\Bar{\R}$ is l.s.c.\ for a.e.\ $t\in[0,T]$, bounded from below on bounded sets, and $t\mapsto\ell(t,x(t),u(t),a(t),\dot x(t),\dot u(t),\dot a(t))$ is summable on $[0,T]$ for each feasible triple $(x(t),u(t),a(t))$.\vspace*{0.02in}

It is proved in \cite[Proposition~2.9]{ve} that assumptions (H2) and (H3) with the constants therein guarantee the $\eta$-uniform prox-regularity of $C(t)$ from \eqref{e:mset} for each $t\in[0,T]$ with $\eta:=\dfrac{\al}{M_3\be}$. We also mention the following well-posedness result established in \cite[Theorem~1]{et}, which says that for any $u(\cdot)\in W^{1,2}([0,T];\R^n)$ and $a(\cdot)\in W^{1,2}([0,T];\R^d)$ under the validity of (H1)--(H3) there exists the unique solution $x(\cdot)\in W^{1,2}([0,T];\R^n)$ to \eqref{e:psp} and \eqref{e:mset} generated by $(u(\cdot),a(\cdot))$. If in addition to the standing assumptions made, the integrand $\ell$ in \eqref{e:Bolza} is convex with respect to the velocity variables $(\dot x,\dot u,\dot a)$ and $\ell(t,\cdot)$ is majorized by a summable function while $\{\dot u^k(\cdot)\}$ is bounded in $L^2([0,T];\R^n)$ and $\{a^k(\cdot)\}$ is bounded in $W^{1,2}([0,T];\R^d)$ along a minimizing sequence of $z^k(\cdot)=(x^k(\cdot),y^k(\cdot),z^k(\cdot))$, then the sweeping control problem $(P)$ admits an {\em optimal solution} in $W^{1,2}([0,T];\R^{2n+d})$; see \cite[Theorem~4.1]{cm3}.\vspace*{-0.05in}

Note that the latter additional assumptions (including the convexity of $\ell$) are not needed for deriving the main results of \cite{cm3} on necessary optimality conditions for problem $(P)$. What actually is needed is a certain local {\em relaxation stability} when a given local minimizer for $(P)$ maintains its local optimality under the convexification of $(P)$ with respect to velocity variables. We are not going to discuss here relaxation procedures of this (Bogolyubov-Young) type while referring the reader to \cite{cm3,dfm,et,m-book2,t} for more details and references. For the purposes of the current paper it is enough mentioning just two results ensuring the relaxation stability of {\em strong} local minimizers (i.e., with respect to a ${\cal C}[0,T]$-neighborhood) in the framework of $(P)$ under our standing assumptions. It holds if either there are only $a$-controls in $(P)$ (see \cite[Theorem~2]{et}), or the set $C$ in \eqref{e:mset} is convex (see \cite[Theorem~4.2]{t}). In fact, in \cite{cm3} we considered a more general type of ``intermediate local minimizers" for $(P)$, but we do not include it here while having in mind efficient applications to the crowd motion model.\vspace*{-0.05in}

Among various results on the study of the optimal control problem $(P)$ obtained in \cite{cm3}, only those related to direct applications to the planar crowd motion model are selected to be presented and discussed in this section. The main ones are on necessary optimality conditions for problem $(P)$ and its slight modification needed for our applications. It means that we put aside the very discrete approximation method of deriving these conditions and the corresponding results on its convergence, necessary optimality conditions for discrete-time problems, and calculating first-order and second-order generalized differential constructions, which are crucial for the method implementation. All of this can be found in \cite{cm3}.\vspace*{-0.05in}

To formulate the necessary optimality conditions for strong local minimizers of $(P)$, we only need to recall the  following (first-order) subdifferential notion for lower semicontinuous (l.s.c.) functions. Given a function $f\colon\R^n\to(-\infty,\infty]$ finite at $\ox$ and l.s.c.\ around this point, the subdifferential of $f$ at $\ox$ is generated by the limiting normal cone \eqref{e:Mor-nc} as
\begin{eqnarray*}
\partial f(\ox):=\big\{v\in\R^n\big|\;(v,-1)\in N\big((\ox,f(\ox)\big);\epi f)\big\},
\end{eqnarray*}
where $\epi f:=\{(x,\al)\in\R^{n+1}|\;\al\ge f(x)\}$. If $f$ is continuous around $\ox$, then we equivalently have
\begin{eqnarray}\label{e:sub}
\partial f(\ox)=\Big\{v\in\R^n\Big|\;\exists\,x_k\to\ox,\;v_k\to v\;\mbox{ with }\;\disp\liminf_{x\to x_k}\frac{f(x)-f(x_k)-\la v_k,x-x_k\ra}{\|x-x_k\|}\ge 0,\; k\in\N\Big\}.
\end{eqnarray}
It is easy to observe from \eqref{e:sub} that the subdifferential mapping $\partial f\colon\R^n\tto\R^n$ is {\em robust}, i.e., its graph
$\gph\partial f:=\{(x,v)\in\R^n\times\R^n|\;v\in\partial f(x)\}$ is closed in $\R^{2n}$. It is not restrictive to assume in what follows that the robustness property keeps holding when the time parameter $t$ is also included in the passage to the limit under subdifferentiation with respect to state and velocity variables in the case of nonautonomous integrands $\ell$ in \eqref{e:Bolza}; see \cite{cm3} for more details. This means the closedness of the set
\begin{eqnarray*}
\big\{(t,z,\dot z,v,w)\big|\;(v,w)\in\partial\ell(t,z,\dot x)\big\},
\end{eqnarray*}
where the subdifferential of $\ell(t,z,\dot z)$, for $z=(x,u,a)$ and $\dot z=(\dot x,\dot u,\dot a)$ is taken with respect to $(z,\dot z)$.\vspace*{-0.05in}

The following theorem, where $N_\O(x):=N(x;\O)$ and
\begin{eqnarray*}
I(x):=\big\{i\in\{1,\ldots,m\}\big|\;g_i(x)=0\big\},
\end{eqnarray*}
is a certain specification of the main result of \cite{cm3} given in  Theorem~8.1 therein.\vspace*{-0.05in}

\begin{theorem}{\bf (necessary optimality conditions for general sweeping control problems.)}\label{Th:NOC} Let $\oz(\cdot)=(\ox(\cdot),\ou(\cdot),\oa(\cdot))\in W^{2,\infty}([0,T])$ be a strong local minimizer for problem $(P)$ under the validity of all the standing assumptions. Suppose also that $\ell$ is continuous in $t$ a.e.\ on $[0,T]$ and is represented as
\begin{eqnarray}\label{separ}
\ell(t,z,\dot z)=\ell_1(t,z,\dot x)+\ell_2(t,\dot u)+\ell_3(t,\dot a),
\end{eqnarray}
where the local Lipschitz constants of $\ell_1(t,\cdot,\cdot)$ and $\ell_3(t,\cdot)$ are essentially bounded on $[0,T]$ and continuous at a.e.\ $t\in[0,T]$ including $t=0$, and where $\ell_2$ is differentiable in $\dot u$ on $\R^n$ with the estimates
\begin{eqnarray*}
\begin{array}{c}
\|\nabla_{\dot{u}}\ell_2(t,\dot{u})\|\le L\|\dot{u}\|\;\mbox{ and }\;\|\nabla_{\dot{u}}\ell_2(t,\dot{u}_1)-\nabla_{\dot{u}}\ell_2(s,\dot{u}_2)\|\le L|t-s|+L\|\dot{u}_1-\dot{u}_2\|
\end{array}
\end{eqnarray*}
holding for all $t,s\in[0,T]$, $\dot a\in\R^d$, and $\dot u,\dot u_1,\dot u_2\in\R^n$ with a uniform constant $L>0$. Moreover, assume that either the set $C$ in \eqref{e:mset} is convex, or $C(t)\equiv C$ therein with $\ell_1=\ell(t,x,a,\dot x)$ and $\ell_2=\ell_2(t)$ in \eqref{separ}. Then there exist dual elements $\lm\ge 0$, $p(\cdot)=(p^x(\cdot),p^u(\cdot),p^a(\cdot))\in W^{1,2}([0,T];\R^n\times\R^n\times\R^d)$, $w(\cdot)=(w^x(\cdot),w^u(\cdot),w^a(\cdot))\in L^2([0,T];\R^{2n+d})$, and $v(\cdot)=(v^x(\cdot),v^u(\cdot),v^a(\cdot))\in L^2([0,T];\R^{2n+d})$ satisfying
\begin{eqnarray}\label{subg}
\big(w(t),v(t)\big)\in\co\partial\ell\big(t,\oz(t),\dot{\oz}(t)\big)\;\mbox{ a.e. }\;t\in[0,T]
\end{eqnarray}
as well as measures $\gamma=(\gg_1,\ldots,\gg_n)\in C^*([0,T];\R^n)$, $\xi^1\in C^*([0,T];\R_+)$, and $\xi^2\in C^*([0,T];\R_-)$ on $[0,T]$ such that the following conditions hold:\\[1ex]
$\bullet$ {\sc Primal-dual dynamic relationships:}
\begin{eqnarray}\label{e:pddr}
\dot{\ox}(t)+f\big(\ox(t),\oa(t)\big)=\sum_{i=1}^m\eta_i(t)\nabla g_i\big(\ox(t)-\ou(t)\big)\;\mbox{ for a.e. }\;t\in[0,T]
\end{eqnarray}
with $\eta(\cdot)\in L^2([0,T];\R^+)$ a.e.\ uniquely determined by representation \eqref{e:pddr} and well defined at $t=T$;
\begin{eqnarray*}
\dot p(t)=\lm w(t)+\big(\nabla_xf(\ox(t),\oa(t))^*(\lm v^x(t)-q^x(t)),0,\nabla_af(\ox(t),\oa(t))^*(\lm v^x(t)-q^x(t))\big),
\end{eqnarray*}
\begin{eqnarray*}
q^u(t)=\lm\nabla_{\dot u}\ell\big(t,\dot{\ou}(t)\big),\;\;q^a(t)\in\lm\partial_{\dot a}\ell_3\big(t,\dot{\oa}(t)\big)\;\mbox{ a.e. }\;t\in[0,T],
\end{eqnarray*}
where $q=(q^x,q^u,q^a)\colon[0,T]\to\R^n\times\R^n\times\R^d$ is a vector function of bounded variation, and its left-continuous representative is given for all $t\in[0,T]$, except at most a countable subset, by
\begin{eqnarray*}
q(t):=p(t)-\int_{[t,T]}\left(-d\gg(s),2\ou(s)d(\xi^1(s)+\xi^2(s))+d\gg(s),0\right).
\end{eqnarray*}
Furthermore, for a.e.\ $t\in[0,T]$ including $t=T$ and for all $i=1,\ldots,m$ we have
\begin{eqnarray}\label{e:implications}
g_i\big(\ox(t)-\ou(t)\big)>0\Lto\eta_i(t)=0,\;\;\eta_i(t)>0\Lto\la\nabla g_i\big(\ox(t)-\ou(t),\lm v^x(t)-q^x(t))\big\ra=0.
\end{eqnarray}
$\bullet$ {\sc Transversality conditions}
\begin{eqnarray*}
\left\{\begin{array}{ll}
-p^x(T)+\disp\sum_{i\in I(\ox(T)-\ou(T))}\eta_i(T)\nabla g_i\big(\ox(T)-\ou(T)\big)\in\lm\partial\ph\big(\ox(T)\big),\quad p^a(T)=0,\\[1ex]
p^u(T)\disp-\sum_{i\in I(\ox(T)-\ou(T))}\eta_i(T)\nabla g_i(\ox(T)-\ou(T))\in
-2\ou(T)\big((N_{[0,r_2]}(\|\ou(T)\|)+N_{[r_1,\infty)}(\|\ou(T)\|)\big)\\[1ex]
\end{array}\right.
\end{eqnarray*}
with the validity of the inclusion
\begin{eqnarray*}
-\sum_{i\in I(\ox(T)-\ou(T))}\eta_i(T)\nabla g_i\big(\ox(T)-\ou(T)\big)\in N_C\big(\ox(T)-\ou(T)\big).
\end{eqnarray*}
$\bullet$ {\sc Measure nonatomicity conditions:}\\
{\bf(a)} Take $t\in[0,T]$ with $g_i(\ox(t)-\ou(t))>0$ whenever $i=1,\ldots,m$. Then there is a neighborhood $V_t$ of $t$ in $[0,T]$ such that $\gg(V)=0$ for all the Borel subsets $V$ of $V_t$.\\
{\bf(b)} Take $t\in[0,T]$ with $r_1<\|\ou(t)\|<r_2$. Then there is a neighborhood $W_t$ of $t$ in $[0,T]$ such that $\xi^1(W)=0$ and $\xi^2(W)=0$ for all  the Borel subsets $W$ of $W_t$.\\[1ex]
$\bullet$ {\sc Nontriviality conditions:} We always have:
\begin{eqnarray*}
\lm+\|q^u(0)\|+\|p(T)\|+\|\xi^1\|_{TV}+\|\xi^2\|_{TV}>0.
\end{eqnarray*}
Furthermore, the following implications hold while ensuring the {\sc enhanced nontriviality}:
\begin{eqnarray*}
\big[g_i\big(x_0-\ou(0)\big)>0,\;i=1,\ldots,m\big]\Longrightarrow\big[\lm+\|p(T)\|+\|\xi^1\|_{TV}+\|\xi^2\|_{TV}>0\big],
\end{eqnarray*}
\begin{eqnarray*}
\big[g_i\big(\ox(T)-\ou(T)\big)>0,\;r_1<\|\ou(T)\|<r_2,\;i=1,\ldots,m\big]\Longrightarrow\big[\lm+\|q^u(0)\|+\|\xi^1\|_{TV}+\|\xi^2\|_{TV}>0\big],
\end{eqnarray*}
where $\|\xi\|_{TV}$ stands for the measure total variation on $[0,T]$.
\end{theorem}\vspace*{-0.05in}

Next we consider a modification of problem $(P)$, where optimization is performed over the pairs $(x(\cdot),a(\cdot))\in W^{1,2}([0,T];\R^{n+d})$
under a fixed $u$-control $\ou(\cdot)\in W^{1,2}([0,T];\R^n)$. The following statement is actually a consequence of Theorem~\ref{Th:NOC} that occurs to be the most appropriate for applications to the controlled planar crowd motion model in the subsequent sections of the paper.\vspace*{-0.05in}

\begin{corollary}{\bf (necessary conditions for sweeping optimal solutions with controlled perturbations).}\label{cor:NOC} Let $(\ox(\cdot),\oa(\cdot))\in W^{2,\infty}([0,T];\R^{n+d})$ be a strong local minimizer for the following problem:
\begin{eqnarray}\label{cost}
\mbox{minimize }\;J[x,a]:=\varphi(x(T))+\int^T_0\ell\big(t,x(t),a(t),\dot{x}(t),\dot{a}(t)\big)dt
\end{eqnarray}
over all the pairs $(x(\cdot),a(\cdot))\in W^{1,2}([0,T];\R^{n+d})$ satisfying the sweeping differential inclusion
\begin{eqnarray*}
 -\dot x(t)\in N\big(x(t)-\ou(t);C\big)+f\big(x(t),a(t)\big)\;\mbox{ a.e. }\;t\in[0,T],\quad x(0):=x_0\in C
\end{eqnarray*}
with the nonconvex set $C$ taken from \eqref{e:mset} and the implicit state constraints
\begin{eqnarray*}
g_i\big(x(t)-\ou(t)\big)\ge 0\;\mbox{ for all }\;t\in[0,T]\;\mbox{ and }\;i=1,\ldots,m.
\end{eqnarray*}
under the corresponding assumptions of Theorem~{\rm\ref{Th:NOC}}. Then there exist a number $\lm\ge0$, subgradient functions $w(\cdot)=(w^x(\cdot),w^a(\cdot))\in L^2([0,T];\R^{n+d})$ and $v(\cdot)=(v^x(\cdot),v^a(\cdot))\in L^2([0,T];\R^{n+d})$ satisfying \eqref{subg}, an adjoint arc $p(\cdot)=(p^x(\cdot),p^a(\cdot))\in W^{1,2}([0,T];\R^{n+d})$, and a measure $\gg=(\gg_1,\ldots,\gg_n)\in C^*([0,T];\R^n)$ on $[0,T]$ such that we have conditions \eqref{e:pddr} and \eqref{e:implications} with the functions $\eta_i(\cdot)\in L^2([0,T];\R_+)$ uniquely defined by representation \eqref{e:pddr} together with the following relationships holding for a.e.\ $t\in[0,T]$:
\begin{eqnarray*}
\left\{\begin{array}{ll}
\dot p^x(t)=\lm w^x(t)+\nabla_x f\big(\ox(t),\oa(t)\big)^*\big(\lm v^x(t)-q^x(t)\big),\\
\dot p^a(t)=\lm w^a(t)+\nabla_a f\big(\ox(t),\oa(t)\big)^*\big(\lm v^x(t)-q^x(t)\big),
\end{array}\right.
\end{eqnarray*}
where the vector function $q(\cdot)=(q^x(\cdot),q^a(\cdot))\colon[0,T]\to\R^n\times\R^d$ is of bounded variation on $[0,T]$ satisfying
\begin{eqnarray*}
q^a(t)\in\lm\partial\ell_3\big(t,\dot{\oa}(t)\big)\;\mbox{ for a.e. }\;t\in[0,T]\;\mbox{ and}
\end{eqnarray*}
\begin{eqnarray*}
q(t):=p(t)+\int_{[t,T]}(d\gg(s),0)
\end{eqnarray*}
with the latter equality holding for the function $q(\cdot)$ of bounded variation everywhere on $[0,T]$ except at most a countable subset for its left-continuous representative. We also have the measure nonatomicity condition {\rm (a)} of Theorem~{\rm\ref{Th:NOC}} and the right-end transversality relationships
\begin{eqnarray*}
\disp-p^x(T)+\sum_{i\in I(\ox(T)-\ou(T))}\eta_i(T)\nabla g_i\big(\ox(T)-\ou(T)\big)\in\lm\partial\ph\big(\ox(T)\big),\quad p^a(T)=0
\end{eqnarray*}
with the validity of the refined nontriviality condition $\lm+\|p(T)\|>0$.
\end{corollary}\vspace*{-0.25in}

\section{Optimization of the Planar Crowd Motion Model}
\setcounter{equation}{0}\vspace*{-0.1in}

In this section we start applications of the obtained results on optimal control of the nonconvex sweeping process presented above to the well-recognized {\em crowd motion model} on the {\em plane}, which is more realistic in practice and much more challenging mathematically in comparison with the corridor version studied in our previous publication \cite{cm2}. In our description of the crowd motion dynamics we follow Maury and Venel \cite{mv2,ve} who developed a mathematical framework for an uncontrolled microscopic model of this type in the form of a sweeping process and provided its numerical simulations with various applications. Our results given in Section~2 allow us to study optimal control of the planar crowd motion model with establishing verifiable conditions for its solution in the general case of $n\ge 2$ participants and their complete realization in the case of $n=2$. We also present several examples of the usage of the obtained results in the characteristic situations of participant interactions on the plane. This section is devoted to the general setting of the controlled crowd motion model with finitely many participants on the plane, while Section~2 deals with solving the formulated optimal control problems for it with two participants.\vspace*{-0.02in}

The {\em microscopic} version of the crowd motion model is based on the following {\em two principles}. On the one hand, each individual has a {\em spontaneous} velocity that he/she would like to have in the absence of other participants. On the other hand, the {\em actual} velocity must take into account. The latter one is incorporated via a projection of the spontaneous velocity into the set of {\em admissible/feasible} velocities, i.e., those which do not violate certain {\em nonoverlapping} constraints.\vspace*{-0.02in}

In what follows we consider $n$ participants $(n\ge 2)$ in the crowd motion model identified with rigid disks on the plane $\R^2$ of the same radius $R$. The center of the $i$-th disk is denoted by $x_i\in\R^2$. Since overlapping is forbidden, the vector of possible positions $x=(x_1,\ldots,x_n)\in\R^{2n}$ has to belong to the set of {\em feasible configurations} defined by
\begin{eqnarray}\label{e:feasible con}
Q:=\big\{x\in\R^{2n}\big|\;D_{ij}(x)\ge 0\;\mbox{ for all }\;i\not=j\big\},
\end{eqnarray}
where $D_{ij}(x):=\|x_i-x_j\|-2R$ is the signed distance between disks $i$ and $j$. Assuming that the participants exhibit the same behavior, their {\em  spontaneous velocities} can be written as
\begin{eqnarray*}
U(x):=\big(U_0(x_1),\ldots,U_0(x_n)\big)\in\R^{2n}\;\mbox{ as }\;x\in Q,
\end{eqnarray*}
where $Q$ is taken from \eqref{e:feasible con}. Observe that the {\em nonoverlapping condition} in \eqref{e:feasible con} does not allow the participants to move with their spontaneous velocity, and the distance between two participants in contact can only increase. To reflect this situation, we introduce the set of {\em feasible velocities} defined by
\begin{eqnarray*}
C_x:=\big\{v=(v_1,\ldots,v_n)\in\R^{2n}\big|\;D_{ij}(x)=0\Lto\big\la G_{ij}(x),v\big\ra\ge 0\;\mbox{ for all }\;i<j\big\},
\end{eqnarray*}
where we use the calculation and notation
\begin{eqnarray*}
G_{ij}(x):=\nabla D_{ij}(x)=\big(0,\ldots,0,-e_{ij}(x),0,\ldots,0,e_{ij}(x),0,\ldots,0\big)\in\R^{2n},\quad\disp e_{ij}(x):=\dfrac{x_j-x_i}{\|x_j-x_i\|}.
\end{eqnarray*}
The {\em actual velocity field} is defined via the Euclidean projection of the spontaneous velocity $U(x)$ at the position $x$ into the feasible velocity set $C_x$ by
\begin{eqnarray*}
\dot x(t)=\Pi\big(U(x);C_x\big)\;\mbox{ for a.e. }\;t\in[0,T]\;\mbox{ with }\;x(0)=x_0\in Q,
\end{eqnarray*}
where $T>0$ is a fixed duration of the process, and where $x_0$ indicates the staring position of the participants. Using the orthogonal decomposition via the sum of mutually polar cone as in \cite{ve}, we get
\begin{eqnarray*}
U(x)\in N_x+\dot{x}(t)\;\mbox{ for a.e. }\;t\in[0,T],\quad x(0)=x_0,
\end{eqnarray*}
where $N_x:=N(x;Q)$ stands for the normal cone to $Q$ at $x$ and can be described in this case as the polar to the feasible velocity set $C_x$ as follows:
\begin{eqnarray*}
N_x=C^*_x=\big\{w\in\R^n\big|\;\la w,v\ra\le 0\;\mbox{ for all }\;v\in C_x\big\},\quad x\in Q.
\end{eqnarray*}
\vspace*{-0.22in}

Let us now rewrite this model in the form used in the version of problem $(P)$ considered in Corollary~\ref{cor:NOC}. To proceed, we specify the nonpolyhedral set $C$ in \eqref{e:mset} by
\begin{eqnarray}\label{e:specified C}
C:=\big\{x=(x_1,\ldots,x_n)\in\R^{2n}\big|\;g_{ij}(x)\ge 0\;\mbox{ for all }\;i\not=j\;\mbox{ as }\;i,j=1,\ldots,n\big\}
\end{eqnarray}
with the signed distance functions $g_{ij}(x):=D_{ij}(x)=\|x_i-x_j\|-2R$. Assume in this framework that all the participants exhibit the same behavior and want to reach the exit by the shortest path. Then their spontaneous velocities can be represented as
\begin{eqnarray*}
U(x)=\big(U_0(x_1),\ldots,U_0(x_n)\big)\;\;\mbox{with}\;\;U_0(x):=-s\nabla D(x)
\end{eqnarray*}
where $D(x_i):=\|x_i\|$ stands for the distance between the position $x_i$ and the exit, and where the scalar $s\ge 0$ denotes the speed. Due to $x\not=0$, we get $\|\nabla D(x)\|=1$ and hence $s=\|U_0(x)\|$. Using it and remembering that each participant tends to maintain his/her spontaneous velocity until reaching the exit, the original perturbation force in this model is described by
\begin{eqnarray*}
f(x):=\big(-s_1\cos\th_1,-s_1\sin\th_1,\ldots,-s_n\cos\th_n,-s_n\sin\th_n\big)\in\R^{2n}\;\mbox{ as }\;x=(x_1,\ldots,x_n)\in Q,
\end{eqnarray*}
where $s_i$ indicates the speed of the participant and $\th_i$ denotes the {\em direction} (i.e., the {\em smallest positive angle} in standard position formed by the positive $x$-axis) of participant $i$ as $i=1,\ldots,n$; see Figure~1.

\begin{center}
\begin{tikzpicture}
\draw [thick, red](-5,0)--(5,0);
\draw [thick, red] (0,0)--(0,6);
\draw [ultra thick, fill = orange] (-0.5,1) rectangle (0.5,0);
\draw node[below] at (0,0){Exit};
\draw node[below] at (0,-1){{\bf Figure 1}};
\draw [thick] (-4,4) circle (0.5);
\draw node[below] at (-4,4){$x_i$};
\draw [ultra thick, blue] [->] (-4,4) -- (-3.2,3.2);
\coordinate[label=above left:] (A) at (-4,4);
\coordinate[label=above left:] (B) at (-5,5);
\coordinate[label=above left:] (C) at (-2.5,4);
\draw[->,purple] (A)--(B);
\draw[->,purple] (A)--(C);
\markangle{B}{A}{C}{$\th_i$}{2}
\end{tikzpicture}
\end{center}
However, if participants $i$ and $j$ are {\em in contact} in the sense that $\|x_i-x_j\|=2R$, i.e., $g_{ij}(x)=0$, then both of them tend to adjust their velocities in order to maintain the distance at least $2R$ with the one in contact. To regulate the actual velocity of all the participants in the presence of the nonoverlapping condition \eqref{e:feasible con}, we involve {\em control functions} $a(\cdot)=(a_1(\cdot),\ldots,a_n(\cdot))$ into perturbations as follows:
\begin{eqnarray}\label{e:controlled pert}
f\big(x(t),a(t)\big):=\big(s_1a_1(t)\cos\th_1(t),s_1a_1(t)\sin\th_1(t),\ldots,s_na_n(t)\cos\th_n(t),s_na_n(t)\sin\th_n(t)\big)
\end{eqnarray}
on the time interval $[0,T]$. To represent this model in the  sweeping control form of Corollary~\ref{cor:NOC}, define recurrently the vector function $\ou=(\ou_1,\ldots,\ou_n)\colon[0,T]\to\R^{2n}$, time independent in our case, by
\begin{equation}\label{e:specified u}
\ou_{i+1}(t)=\ou_i(t):=\left(\dfrac{r}{\sqrt{2n}},\dfrac{r}{\sqrt{2n}}\right),\quad i=1,\ldots,n-1,
\end{equation}
where the number $r>0$ is any constant in the interval $[r_1,r_2]$. It follows from \eqref{e:specified u} that the nonoverlapping condition \eqref{e:feasible con} can be written in the {\em state constraint} form
\begin{eqnarray}\label{e:state constraints}
x(t)-\ou(t)\in C\;\;\mbox{for all}\;\;t\in[0,T],
\end{eqnarray}
where $C$ is taken from \eqref{e:specified C}, and where the points $t\in[0,T]$ with $\|x_j(t)-x_i(t)\|=2R$ are exactly those at which the motion $x(t)-\ou(t)$ hits the nonpolyhedral constraint set \eqref{e:specified C}.\vspace*{-0.05in}

Summarizing the above discussions, we can represent the planar crowd model dynamics in the following form of the constrained controlled sweeping process:
\begin{equation}\label{e:dyn sys}
\left\{\begin{array}{ll}
-\dot x(t)\in N\big(x(t);C(t))+f(x(t),a(t)\big)\;\mbox{ for a.e. }\;t\in[0,T],\\
C(t):=C+\ou(t),\;\|\ou(t)\|=r\in[r_1,r_2]\;\mbox{ on }\;[0,T],\;x(0)=x_0\in C(0),
\end{array}\right.
\end{equation}
with $C,f$, and $\ou(t)$ taken from \eqref{e:specified C}--\eqref{e:specified u} with the state constraints \eqref{e:state constraints} implicitly contained in \eqref{e:dyn sys}.\vspace*{-0.05in}

Now we have a possibility to optimize the controlled crowd motion dynamics \eqref{e:dyn sys} by choosing an appropriate {\em cost functional} of type \eqref{cost}. It seems naturally to associate with \eqref{e:dyn sys} the following cost
\begin{eqnarray}\label{e:cost functional}
\mbox{minimize}\;\;J[x,a]:=\dfrac{1}{2}\left(\|x(T)\|^2+\int^T_0\|a(t)\|^2dt\right),
\end{eqnarray}
which reflects the simultaneous minimization of the {\em distance} of all the participants to the exit at the {\em origin} (see Figure~2) together with the {\em energy} of feasible controls $a(\cdot)$ used to adjust spontaneous velocities.

\begin{center}
\begin{tikzpicture}
\draw [thick, red](-5,0)--(5,0);
\draw [thick, red] (0,0)--(0,6);
\draw [ultra thick, fill = orange] (-0.5,1) rectangle (0.5,0);
\draw node[below] at (0,0){Exit};
\draw [thick] (-4,4) circle (0.5);
\draw node[above] at (-4,4){$x_i$};
\draw [thick] (-3,4) circle (0.5);
\draw node[above] at (-3,4){$x_j$};
\draw [thick] (-2,3) circle (0.5);
\draw node[above] at (-2,3){$x_1$};
\draw [thick] (2,2) circle (0.5);
\draw node[above] at (2,2){$x_n$};
\draw [thick] (4,5) circle (0.5);
\draw node[above] at (4,5){$\ldots$};
\draw [thick] (0,3) circle (0.5);
\draw node[above] at (0,3){$x_{n-1}$};
\draw [ultra thick, blue] [->] (-4,4) -- (-3.5,3.5);
\draw [ultra thick, blue] [->] (-3,4) -- (-2.7,3.6);
\draw [ultra thick, blue] [->] (-2,3) -- (-1.8,2.7);
\draw [ultra thick, blue] [->] (2,2) -- (1.5,1.5);
\draw [ultra thick, blue] [->] (4,5) -- (3,3.75);
\draw [ultra thick, blue] [->] (0,3) -- (0,2.4);
\draw node[below] at (0,-1){{\bf Figure 2}};
\end{tikzpicture}
\end{center}
Thus we arrive at a dynamic optimization problem for the planar crowd motion model formalized via optimal control of the nonpolyhedral sweeping process treated in Corollary~\ref{cor:NOC} over feasible pairs $(\ox(\cdot),\oa(\cdot))\in W^{1,2}([0,T];\R^{3n})$. The existence of optimal solutions in \eqref{e:state constraints}--\eqref{e:cost functional} with the specified data from \eqref{e:specified C}--\eqref{e:specified u} follows from our discussions in Section~2 (cf.\ \cite[Theorem~4.1]{cm3}). Now our goal is to apply the necessary conditions of Corollary~\ref{cor:NOC} to the general crowd motion model under consideration.\vspace*{-0.05in}

\begin{theorem}{\bf (necessary conditions for optimal control of the planar crowd motions).}\label{n-crowd} Let $(\ox(\cdot),\oa(\cdot))\in W^{2,\infty}([0,T];\R^{2n})$ be a strong local minimizer for the crowd motion problem \eqref{e:state constraints}--\eqref{e:cost functional} with the data from \eqref{e:specified C}--\eqref{e:specified u}. Then there exist $\lm\ge 0$, $\eta_{ij}(\cdot)\in L^2([0,T];\R_+)$ $(i,j=1,\ldots,n)$ well defined at $t=T$, $w(\cdot)=(w^x(\cdot),w^a(\cdot))\in L^2([0,T];\R^{2n})$, $v(\cdot)=(v^x(\cdot),v^a(\cdot))\in L^2([0,T];\R^{3n})$, $p(\cdot)=(p^x(\cdot),p^a(\cdot))\in W^{1,2}([0,T];\R^{2n})$, a measure $\gg\in C^*([0,T];\R^{2n})$ on $[0,T]$, and a vector function $q(\cdot)=(q^x(\cdot),q^a(\cdot))\colon[0,T]\to\R^{3n}$ of bounded variation on $[0,T]$ such that the following conditions are satisfied:
\begin{enumerate}
\item[{\bf(1)}] $w(t)=\big(0,\oa(t)\big)$,\quad $v(t)=(0,0)$ for a.e.\ $t\in[0,T]$;
\item[{\bf(2)}]
$\dot\ox(t)+\big(s_1\oa_1(t)\cos\th_1(t),s_1\oa_2(t)\sin\th_2(t),\ldots,s_n\oa_n(t)\cos\th_n(t),s_n\oa_n(t)\sin\th_n(t)\big)$\\[1ex]
$=\disp\sum_{i<j}\eta_{ij}(t)\nabla g_{ij}(\ox(t)-\ou(t))$\\[1ex]
$=\disp\bigg(-\sum_{j>1}\eta_{1j}(t)\dfrac{\ox_j(t)-\ox_1(t)}{\|\ox_j(t)-\ox_1(t)\|},\ldots,\sum_{i<j}\eta_{ij}(t)
\dfrac{\ox_j(t)-\ox_i(t)}{\|\ox_j(t)-\ox_i(t)\|}-\sum_{i>j}\eta_{ji}(t)\dfrac{\ox_i(t)-\ox_j(t)}{\|\ox_i(t)-\ox_j(t)\|},$\\[1ex]
$\disp\ldots,\sum_{j<n}\eta_{jn}(t)\dfrac{\ox_n(t)-\ox_j(t)}{\|\ox_n(t)-\ox_j(t)\|}\bigg)$.
\item[{\bf(3)}] $\|\ox_i(t)-\ox_j(t)\|>2R\Longrightarrow\eta_{ij}(t)=0$ for all $i<j$ and a.e.\ $t\in[0,T]$;
\item[\bf{(4)}] $\eta_{ij}(t)>0\Longrightarrow\left\la q^x_j(t)-q^x_i(t),\ox_j(t)-\ox_i(t)\right\ra=0$ for all $i<j$ and a.e.\ $t\in[0,T]$;
\item[{\bf(5)}]
$\left\{\begin{array}{ll}
\dot p(t)=\bigg(0,\lm\oa_1(t)-s_1\big[\cos\th_1(t)q^x_{11}(t)+\sin\th_1(t)q^x_{12}(t)\big],\\\ldots,\lm\oa_n(t)-s_n\big[\cos\th_n(t)q^x_{n1}(t)+\sin\th_n(t)q^x_{n2}(t)\big]
\bigg)\end{array}\right.$
\item[{\bf(6)}] $q^x(t)=p^x(t)+\gg([t,T])$ for a.e.\ $t\in[0,T]$;
\item[{\bf(7)}] $q^a(t)=p^a(t)=0$ for a.e.\ $t\in[0,T]$;
\item[{\bf(8)}]$\left\{\begin{array}{ll}
p^x(T)+\lm\ox(T)=\disp\bigg(-\sum_{j>1}\eta_{1j}(T)\dfrac{\ox_j(T)-\ox_1(T)}{\|\ox_j(T)-\ox_1(T)\|},\ldots,\\[3ex]
\disp\sum_{i<j}\eta_{ij}(T)\dfrac{\ox_j(T)-\ox_i(T)}{\|\ox_j(T)-\ox_i(T)\|}-\sum_{i>j}\eta_{ji}(T)\dfrac{\ox_i(T)-\ox_j(T)}{\|\ox_i(T)-\ox_j(T)\|},\\[3ex]
\disp\ldots,\sum_{j<n}\eta_{jn}(T)\dfrac{\ox_n(T)-\ox_j(T)}{\|\ox_n(T)-\ox_j(T)\|}\bigg);
\end{array}\right.
$\item[\bf{(9)}] $p^a(T)=0$;
\item[{\bf(10)}] $\lm+\|p^x(T)\|>0$.
\end{enumerate}
\end{theorem}\vspace*{-0.05in}
{\bf Proof.} First we need to make sure that all the assumptions imposed in Corollary~\ref{cor:NOC} hold for problem \eqref{e:state constraints}--\eqref{e:cost functional}. This means checking the standing assumptions (H1)--(H4) together with the additional assumptions formulated in Theorem~\ref{Th:NOC}. Since the other assumptions obviously hold in the framework of Corollary~\ref{cor:NOC}, it remains to verify of conditions \eqref{e:dist-inequality}--\eqref{e:w-inverse-triangle-in} in (H2) and (H3). Indeed, we have that each function $g_{ij}$ is convex, belongs to the space ${\cal C}^2(V_{ij})$ on the open set
\begin{eqnarray*}
V_{ij}:=\big\{x\in\R^{2n}\big|\;\|x_i-x_j\|-R>0\big\}\;\mbox{ for all }\;i,j=1,\ldots,n,
\end{eqnarray*}
and satisfies estimate \eqref{e:dist-inequality} with $c:=\dfrac{R}{\sqrt2}$; see \cite[Proposition~2.9]{ve}. Furthermore, it is obvious that
\begin{eqnarray*}
|\nabla g_{ij}(x)\|=\sqrt2\;\mbox{ and }\;\|\nabla^2g_{ij}(x)\|\le\dfrac{2}{r}\;\mbox{ as }\;x\in V_{ij},
\end{eqnarray*}
and hence the inequalities in \eqref{e:bound-grad} hold. Finally, it follows from \cite[Proposition~4.7]{ve} that there exists $\be>1$, which in fact can be calculated by the formula
\begin{eqnarray*}
\be=\disp3\sqrt2n\left(\dfrac{3}{\sin\left(2\pi/n\right)}\right)^n
\end{eqnarray*}
and ensures the validity of the estimate
\begin{eqnarray*}
\disp\sum_{(i,j)\in I(x)}\al_{ij}\|\nabla g_{ij}(x)\|\le\be\left\|\sum_{(i,j)\in I(x)}\al_{ij}\nabla g_{ij}(x)\right\|\;\mbox{ for all }\;x\in Q
\end{eqnarray*}
with $I(x):=\{(i,j)|\;g_{ij}(x)=0,\;i<j\}$ and $\al_{ij}\ge 0$. This shows that assumption \eqref{e:w-inverse-triangle-in} is satisfied.\vspace*{-0.05in}

Dealing with the given data \eqref{e:specified C}--\eqref{e:specified u} of the crowd motion control; problem \eqref{e:state constraints}--\eqref{e:cost functional}, it is not hard to verify directly that the necessary optimality conditions obtained Corollary~\ref{cor:NOC} reduce in the case under consideration to those presented in items (1)--(10) of the theorem. $\h$

Let us now elaborate the necessary optimality conditions of Theorem~\ref{n-crowd}. As discussed above, when two participants $i$ and $j$ are in contact in the sense that $\|\ox_i(t_{ij})-\ox_j(t_{ij})\|=2R$ for some contact time $t_{ij}\in[0,T]$, they tend to adjust their velocities in order to keep the distance at least $2R$ with the one in contact. Thus it is natural to accept that both participants $i$ and $j$ maintain their new constant velocities after the time $t=t_{ij}$ until either reaching someone or the end of the process at time $t=T$ and that the control functions are constant, i.e., $\oa_i(t)\equiv\oa_i$ as $i=1,\ldots,n$. It tells us that the velocities of all the participants are piecewise constant on $[0,T]$ in this setting.\vspace*{-0.05in}

Further, we can represent the differential equations in (2) as
\begin{eqnarray*}
\left\{\begin{array}{ll}
\dot\ox_1(t)=-s_1\oa_1\big(\cos\th_1(t),\sin\th_1(t)\big)-\disp\sum_{j>1}\eta_{1j}(t)\dfrac{\ox_j(t)-\ox_1(t)}{\|\ox_j(t)-\ox_1(t)\|},\\[1ex]
\dot\ox_i(t)=-s_i\oa_i\big(\cos\th_i(t),\sin\th_i(t)\big)+\disp\sum_{i<j}\eta_{ij}(t)\dfrac{\ox_j(t)-\ox_i(t)}{\|\ox_j(t)-\ox_i(t)\|}-\sum_{i>j}\eta_{ji}(t)
\dfrac{\ox_i(t)-\ox_j(t)}{\|\ox_i(t)-\ox_j(t)\|},\;i=2,\ldots,n-1,\\[1ex]
\dot\ox_n(t)=-s_n\oa_n\big(\cos\th_n(t),\sin\th_n(t)\big)+\disp\sum_{j<n}\eta_{jn}(t)\dfrac{\ox_n(t)-\ox_j(t)}{\|\ox_n(t)-\ox_j(t)\|}
\end{array}\right.
\end{eqnarray*}
for a.e.\ $t\in[0,T]$. It follows from the optimality conditions in (5) and (7) that
\begin{eqnarray}\label{e:a-q}
\lm\oa_i=s_i\big[\cos\th_i(t)q^x_{i1}(t)+\sin\th_i(t)q^x_{i2}(t)\big]\;\mbox{ for a.e. }\;t\in[0,T]\;\mbox{ and all }\;i=1,\ldots,n.
\end{eqnarray}

Fix now $i,j\in\{1,\ldots,n\}$, and let $t_{ij}$ be the first time when $\|\ox_i(t_{ij})-\ox_j(t_{ij})\|=2R$. For such indices $i$ and $j$ define the positive numbers
\begin{eqnarray}\label{e:ct}
\vTh^{ij}:=\min\big\{t_{i'j'}\big|\;t_{i'j'}>t_{ij},\;\nu=1,\ldots,n\big\},\quad\vTh_{ij}:=\max\big\{t_{i'j'}\big|\;t_{i'j'}<t_{ij},\;\nu=1,\ldots,n\big\}
\end{eqnarray}
If $\eta_{ij}(t)>0$ for some $i,j\in\{1,\ldots,n-1\}$ and $t\in[0,T]$, we deduce from (4) that
\begin{eqnarray}\label{e:qij}
\big(q^x_{j1}(t)-q^x_{i1}(t)\big)\big(\ox_{j1}(t)-\ox_{i1}(t)\big)+\big(q^x_{j2}(t)-q^x_{i2}(t)\big)\big(\ox_{j2}(t)-\ox_{i2}(t)\big)=0.
\end{eqnarray}
After the contact time $t=t_{ij}$, the two participants $i$ and $j$ adjust to the same velocity and maintain the new velocity until either reaching someone or at the end of the process. Thus we have $\dot\ox_i(t)=\dot\ox_j(t)$ and  $\|\ox_i(t)-\ox_j(t)\|=2R$ for all $t\in[t_{ij},\vTh^{ij})$. As a consequence, it yields
\begin{eqnarray*}
\begin{aligned}
\ox_i(t)-\ox_j(t)&=\ox_i(0)-\ox_j(0)+\int^t_0\big[\dot\ox_i(s)-\dot\ox_j(s)\big]ds=\ox_i(0)-\ox_j(0)+\int^{t_{ij}}_0\big[\dot\ox_i(s)-\dot\ox_j(s)\big]ds\\
&=\ox_i(t_{ij})-\ox_j(t_{ij})=2R(\cos\th_{ij},\sin\th_{ij})\;\mbox{ for all }\;t\in[t_{ij},\vTh^{ij}),
\end{aligned}
\end{eqnarray*}
where $\th_{ij}$ indicates the direction of the vector $\ox_i(t_{ij})-\ox_j(t_{ij})$; see Figure~3.

\begin{center}
\begin{tikzpicture}
\draw [thick, red](-5,0)--(5,0);
\draw [thick, red] (0,0)--(0,6);
\draw [ultra thick, fill = orange] (-0.5,1) rectangle (0.5,0);
\draw node[below] at (0,0){Exit};
\draw node[below] at (0,-1){{\bf Figure 3}};
\draw [thick] (-4,4) circle (0.5);
\draw node[below] at (-4,4){$x_i$};
\draw [thick] (-3.5,3.13) circle (0.5);
\draw node[below] at (-3.5,3.13){$x_j$};
\coordinate[label=above left:] (A) at (-4,4);
\coordinate[label=above left:] (B) at (-3.5,3.13);
\coordinate[label=above left:] (C) at (-2.5,3.13);
\draw[->,purple] (B)--(A);
\draw[->,purple] (B)--(C);
\markangle{A}{B}{C}{$\th_{ij}$}{2}
\end{tikzpicture}
\end{center}
If $\eta_{ij}(t_{ij})>0$ and if the two participants have the same direction at the contact time, i.e., $\th_i(t_{ij})=\th_j(t_{ij})$, then it follows from \eqref{e:a-q} and \eqref{e:qij} that
\begin{eqnarray*}
\big(q^x_{j1}(t_{ij})-q^x_{i1}(t)\big)2R\cos\th_i(t_{ij})+\big(q^x_{j2}(t)-q^x_{i2}(t_{ij})\big)2R\sin\th_i(t_{ij})=0,
\end{eqnarray*}
which is equivalent to the equation
$$
\cos\th_j(t_{ij})q^x_{j1}(t_{ij})+\sin\th_j(t_{ij})q^x_{j2}(t_{ij})=\cos\th_i(t_{ij})q^x_{i1}(t_{ij})+\sin\th_i(t_{ij})q^x_{i2}(t_{ij}).
$$
Hence we have $\lm s_j\oa_i=\lm s_i\oa_j$ and therefore arrive at the equality
\begin{eqnarray} \label{e:aij}
s_j\oa_i=s_i\oa_j\;\mbox{ for }\;i,j\in\big\{1,\ldots,n\big\}
\end{eqnarray}
provided that $\lm\not=0$; otherwise we do not have enough information to proceed.\vspace*{-0.02in}

To describe the relation between the positions of participants $i$ and $j$ during the period of contact, we introduce the function $d_{ij}(t):=\dfrac{\ox_i(t)-\ox_j(t)}{\|\ox_i(t)-\ox_j(t)\|}$, which admits the following representation:
\begin{equation}\label{e:dij}
d_{ij}(t)=\chi_{[t_{ij},\vTh^{ij})}(t)\big(\cos\th_{ij}(t_{ij}),\sin\th_{ij}(t_{ij})\big),
\end{equation}
where $\vTh^{ij}$ are taken from \eqref{e:ct}, and where $\chi_S(t)$ stands for the characteristic function of the set $S$ that is equal to $1$ if $t\in S$ and to $0$ otherwise.\vspace*{-0.02in}

With the usage of \eqref{e:dij} and the discussions above, the crowd motion differential equations from (2) can be written in the form
\begin{eqnarray}
\label{e:velocities1}
\left\{\begin{array}{ll}
\dot\ox_1(t)=-s_1\oa_1\big(\cos\th_1(t),\sin\th_1(t)\big)-\disp\sum_{j>1}\eta_{1j}(t)d_{j1}(t),\\[1ex]
\dot\ox_i(t)=-s_i\oa_i\big(\cos\th_i(t),\sin\th_i(t)\big)+\disp\sum_{i<j}\eta_{ij}(t)d_{ji}(t)-\sum_{i>j}\eta_{ji}(t)d_{ij}(t),\;\;i=2,\ldots,n-1,\\[1ex]
\dot\ox_n(t)=-s_n\oa_n\big(\cos\th_n(t),\sin\th_n(t)\big)+\disp\sum_{j<n}\eta_{jn}(t)d_{nj}(t).
\end{array}\right.
\end{eqnarray}\vspace*{-0.1in}

In the next section we consider the crowd motion control problem for two participants and show how to explicitly solve the problem using the necessary optimality conditions obtained in theory.\vspace*{-0.15in}

\section{Crowd Motion Control Problem with Two Participants}
\setcounter{equation}{0}\vspace*{-0.1in}

We first proceed with deriving analytic relations for optimal solutions to the crowd motion control problem \eqref{e:state constraints}--\eqref{e:cost functional} in the case of $n=2$ participants that allow us to completely solve the problem under consideration. Then we provide numerical calculations in the most characteristic settings reflecting possible interactions between the model participants.\vspace*{-0.02in}

It follows from \eqref{e:velocities1} that the velocities of two participants before and after the contact time $t_{12}$ are given, respectively, by the equations
\begin{eqnarray*}
\left\{\begin{array}{ll}
\dot\ox_1(t)=\big(-s_1\oa_1\cos\th_1(0),-s_1\oa_1\sin\th_1(0)\big),\\[1ex]
\dot\ox_2(t)=\big(-s_2\oa_2\cos\th_2(0),-s_2\oa_2\sin\th_2(0)\big)
\end{array}\right.
\;\mbox{ for }\;t\in[0,t_{12});
\end{eqnarray*}
\begin{eqnarray*}
\left\{\begin{array}{ll}
\dot\ox_1(t)=-s_1\oa_1\big(\cos\th_1(0),\sin\th_1(0)\big)-\eta_{12}(t)d_{21}(t),\\[1ex]
\dot\ox_2(t)=-s_2\oa_2\big(\cos\th_2(0),\sin\th_2(0)\big)+\eta_{12}(t)d_{21}(t)
\end{array}\right.
\;\mbox{ for }\;t\in[t_{12},T].
\end{eqnarray*}
Since $d_{21}(t)=(\cos\th_{21},\sin\th_{21})$ for $t\in[t_{12},T]$, the function $\eta_{12}(\cdot)$ is calculated by
\begin{eqnarray*}
\eta_{12}(t)=\left\{\begin{array}{ll}
\eta_{12}(0)=0,\;\;&t\in[0,t_{12}),\\[1ex]
\eta_{12}(t_{12})=:\eta_{12},\;\;&t\in[t_{12},T],
\end{array}\right.
\end{eqnarray*}
on $[0,T]$, and hence we arrive at the velocity representations on $[t_2,T]$:
\begin{eqnarray*}
\left\{\begin{array}{ll}
\dot\ox_1(t)=\big(-s_1\oa_1\cos\th_1(0)-\eta_{12}\cos\th_{21},-s_1\oa_1\sin\th_1(0)-\eta_{12}\sin\th_{21}\big),\\[1ex]
\dot\ox_2(t)=\big(-s_2\oa_2\cos\th_2(0)+\eta_{12}\cos\th_{21},-s_2\oa_2\sin\th_2(0)+\eta_{12}\sin\th_{21}\big).
\end{array}\right.
\end{eqnarray*}
Since the two participants have the same velocities at $t=t_{12}$ and keep the new velocity until the end of the process, it follows that $\dot\ox_1(t)=\dot\ox_2(t)$ for all $t\in[t_{12},T]$, which yields
\begin{eqnarray}\label{e:eta cal}
\begin{cases}
(\cos\th_{21})\eta_{12}=\dfrac{s_2\oa_2\cos\th_2(0)-s_1\oa_1\cos\th_1(0)}{2},\\[1ex]
(\sin\th_{21})\eta_{12}=\dfrac{s_2\oa_2\sin\th_2(0)-s_1\oa_1\sin\th_1(0)}{2}.
\end{cases}
\end{eqnarray}
This allows us to calculate the value of $\eta_{12}$ by
\begin{equation}
\label{e:eta cal1}
\eta_{12}=\frac{1}{2}\sqrt{\big(s_1^2\oa_1^2+s_2^2\oa_2^2-2s_1s_2\oa_1\oa_2\cos(\th_1(0)-\th_2(0))\big)}.
\end{equation}
Furthermore, in this way we calculate the function $d_{21}(t)$ for $t\in[t_{12},T]$ by
\begin{eqnarray*}
\begin{aligned}
d_{21}(t)&=\dfrac{\ox_2(t)-\ox_1(t)}{\|\ox_2(t)-\ox_1(t)\|}=\dfrac{1}{2R}\left(\ox_2(0)-\ox_1(0)+\int^t_0(\dot\ox_2(s)-\dot\ox_1(s))ds\right)\\
&=\dfrac{1}{2R}\Big(\ox_{21}(0)-\ox_{11}(0)-2t_{12}\cos\th_{21}\eta_{12},\ox_{22}(0)-\ox_{12}(0)-2t_{12}\sin\th_{21}\eta_{12}\Big).
\end{aligned}
\end{eqnarray*}
On the other hand, we have $d_{21}(t_{12})=(\cos\th_{21},\sin\th_{21})$ from \eqref{e:dij}, and so
\begin{eqnarray*}
\begin{cases}
\cos\th_{21}=\dfrac{\ox_{21}(0)-\ox_{11}(0)}{2t_{12}\eta_{12}+2R},\\[2ex]
\sin\th_{21}=\dfrac{\ox_{22}(0)-\ox_{12}(0)}{2t_{12}\eta_{12}+2R}.
\end{cases}
\end{eqnarray*}
This immediately implies the relationship
\begin{eqnarray}\label{e:t12}
t_{12}\eta_{12}=\dfrac{\|\ox_2(0)-\ox_1(0)\|-2R}{2},
\end{eqnarray}
which leads us to the calculation of $\cos\th_{21}$ and $\sin\th_{21}$ as follows:
\begin{eqnarray}\label{e:dir21}
\cos\th_{21}=\dfrac{\ox_{21}(0)-\ox_{11}(0)}{\|\ox_1(0)-\ox_2(0)\|},\quad\sin\th_{21}=\dfrac{\ox_{22}(0)-\ox_{12}(0)}{\|\ox_1(0)-\ox_2(0)\|}.
\end{eqnarray}

Now we can summarize the above discussions and present computation formulas for determining optimal trajectories of the motions depending on the parameters involved.\vspace*{-0.05in}

\begin{proposition}{\bf (calculation of optimal trajectories in motions with two participants).}\label{2-crowd} Optimal trajectories in the controlled motion problem \eqref{e:state constraints}--\eqref{e:cost functional} with $n=2$ are determined by
\begin{eqnarray*}
\left\{\begin{array}{ll}
\ox_1(t)=\big(\ox_{11}(0),\ox_{12}(0)\big)+\big(-s_1\oa_1t\cos\th_1(0),-s_1\oa_1t\sin\th_1(0)\big),\\[2ex]
\ox_2(t)=\big(\ox_{21}(0),\ox_{22}(0)\big)+\big(-s_2\oa_2t\cos\th_2(0),-s_2t\oa_2\sin\th_2(0)\big)
\end{array}\right.\;\mbox{ for }\;t\in[0,t_{12});
\end{eqnarray*}
\begin{eqnarray*}
\left\{\begin{array}{ll}
\ox_1(t)=\big(\ox_{11}(0),\ox_{12}(0)\big)+\big(-s_1\oa_1t(\cos\th_1(0)-\eta_{12}(t-t_{12})\cos\th_{21},\\
-s_1\oa_1t\sin\th_1(0)-\eta_{12}t-t_{12})\sin\th_{21}\big),\\[2ex]
\ox_2(t)=(\ox_{21}(0),\ox_{22}(0))+\big(-s_2\oa_2t\cos\th_2(0)+\eta_{12}(t-t_{12})\cos\th_{21},\\
-s_2\oa_2t\sin\th_2(0)+\eta_{12}(t-t_{12})\sin\th_{21}\big)
\end{array}\right.\;\mbox{ for }\;t\in[t_{12},T],
\end{eqnarray*}
where $\eta_{12}\ge 0$ is calculated in \eqref{e:eta cal1}, where the contact time $t_{12}$ is taken from  \eqref{e:t12}, where $\sin\th_{21}$ and $\cos\th_{21}$ the directional angle $\th_{12}$ at the contact time can be found from $d_{21}(t_{12})=(\cos\th_{21},\sin\th_{21})$ via \eqref{e:dij}, where $\th_1(0)$ and $\th_2(0)$ are angles of the participant direction at the initial positions, and where the piecewise constraint optimal controls $\oa_1,\oa_2$ with switching at the contact time satisfy \eqref{e:aij} if $\eta_{12}>0$.
\end{proposition}\vspace*{-0.05in}
{\bf Proof.} Follows from the above by integration via the Newton-Leibniz formula. $\h$

Note that in the remaining case of $\eta_{12}=0$ in Proposition~\ref{2-crowd}, we get from \eqref{e:eta cal1} that
\begin{eqnarray}\label{eta0}
s_1^2\oa_1^2+s_2^2\oa_2^2=2s_1s_2\oa_1\oa_2\cos\big(\th_1(0)-\th_2(0)\big),
\end{eqnarray}
which tells us by \eqref{e:t12} that $\|\ox_2(0)-\ox_1(0)\|=2R$. The latter means that the two participants must be in contact at the initial time. Let us analyze this situation on the following numerical example.\vspace*{-0.1in}

\begin{example}{\bf (participants are in contact at the initial time).}\label{ex1} {\rm Consider the optimal control problem in \eqref{e:state constraints}--\eqref{e:cost functional} with the initial data:
\begin{eqnarray*}
n=2,\;T=6,\;s_1=6,\;s_2=3,\;x_{01}= \left(-48-\frac{6}{\sqrt2},48+\dfrac{6}{\sqrt2}\right),\;x_{02}=(-48,48),\;R=3.
\end{eqnarray*}
In this setting we have $t_{12}=0$, $\th_1(0)=\th_2(0)=135^\circ$; see Figure~4.

\begin{center}
\begin{tikzpicture}
\draw [thick, red](-5,0)--(5,0);
\draw [thick, red] (0,0)--(0,6);
\draw [ultra thick, fill = orange] (-0.5,1) rectangle (0.5,0);
\draw node[below] at (0,0){Exit};
\draw node[below] at (0,-1){{\bf Figure 4}};
\draw [thick] (-4,4) circle (0.5);
\draw node[below] at (-4,4){$x_1$};
\draw [thick] (-3.29,3.29) circle (0.5);
\draw node[above] at (-3.29,3.29){$x_2$};
\draw [ultra thick, blue] [->] (-4,4) -- (-3.2,3.2);
\draw [ultra thick, blue] [->] (-3.29,3.29) -- (-2.49,2.49);
\coordinate[label=above left:] (A) at (-4,4);
\coordinate[label=above left:] (B) at (-5,5);
\coordinate[label=above left:] (C) at (-2.5,4);
\draw[->,purple] (A)--(B);
\draw[->,purple] (A)--(C);
\markangle{B}{A}{C}{$135^\circ$}{2}
\end{tikzpicture}
\end{center}
Then $\cos\th_{21}=\frac{\sqrt2}{2}$ and $\sin\th_{21}=-\frac{\sqrt2}{2}$, and thus the optimal trajectories $\ox_1(\cdot)$ and $\ox_2(\cdot)$ are given due to Proposition~\ref{2-crowd} by the formulas
\begin{eqnarray*}
\left\{\begin{array}{ll}
\ox_1(t)=\Big(-48-\frac{6}{\sqrt2}+\big(3\sqrt2\oa_1-\frac{\sqrt2}{2}\eta_{12}\big)t,48+\frac{6}{\sqrt2}+\big(-3\sqrt2\oa_1+\frac{\sqrt2}{2}\eta_{12}\big)t\Big),\\[1ex]
\ox_2(t)=\Big(-48+\big(\frac{3\sqrt2}{2}\oa_2+\frac{\sqrt2}{2}\eta_{12}\big)t,48+\big(-\frac{3\sqrt2}{2}\oa_2-\frac{\sqrt2}{2}\eta_{12}\big)t\Big).
\end{array}\right.
\end{eqnarray*}
There are the following two possible cases to examine:\\[1ex]
{\bf Case~1:}  $\eta_{12}(t_{12})=\eta_{12}(0)=0$. In this case we get from \eqref{eta0} that $\oa_2=2\oa_1$, and so
\begin{eqnarray*}
\begin{cases}
\ox_1(t)=\left(-48-\frac{6}{\sqrt2}+3\sqrt2\oa_1t,48+\frac{6}{\sqrt2}-3\sqrt2\oa_1t\right),\\[1ex]
\ox_2(t)=\left(-48+3\sqrt2\oa_1t,48-3\sqrt2\oa_1t\right).
\end{cases}
\end{eqnarray*}
Substituting it into the cost functional \eqref{e:cost functional} gives us
\begin{eqnarray*}
J[x,a]=1311\oa^2_1-36(96\sqrt2+6)\oa_1+\left(48+\frac{6}{\sqrt2}\right)^2+48^2.
\end{eqnarray*}
It is easy to see that $J$ attains its minimum at $\oa_1=\dfrac{(96\sqrt2+6)18}{1311}\approx 1.95$. Thus $\oa_2=2\oa_1\approx 3.9$, the minimal cost value in this case is $J\approx 66.49$, and the optimal trajectory is calculate on $[0,6]$ by
\begin{eqnarray}\label{traj1}
\ox_1(t)=\left(-48-\frac{6}{\sqrt2}+8.27t,48+\frac{6}{\sqrt2}-8.27t \right),\;\ox_2(t)=\left(-48+8.27t,48-8.27t\right).
\end{eqnarray}
{\bf Case~2:} $\eta_{12}(t_{12})=\eta_{12}(0)>0$. In this case we get $\oa_1=2\oa_2$ from \eqref{e:aij}, and hence $\eta_{12}(t_{12})=\frac{9}{2}\oa_2$. Then the trajectory calculations of Proposition~\ref{2-crowd} read now as
\begin{eqnarray*}
\left\{\begin{array}{ll}
\ox_1(t)=\left(-48-\frac{6}{\sqrt2}+\frac{15\sqrt2}{4}\oa_2t,48+\frac{6}{\sqrt2}-\frac{15\sqrt2}{4}\oa_2t\right),\\[1ex]
\ox_2(t)=\left(-48+\frac{15\sqrt2}{4}\oa_2t,48-\frac{15\sqrt2}{4}\oa_2t\right),
\end{array}\right.
\end{eqnarray*}
and the cost functional is calculated by
\begin{eqnarray*}
J[x,a]=2040\oa^2_2-45(96\sqrt2+6)\oa_2+\left(48+\frac{6}{\sqrt2}\right)^2+48^2
\end{eqnarray*}
It attains its minimum at $\oa_2=\dfrac{45(96\sqrt2+6)}{4080}\approx1.56$, and hence $\oa_1=2\oa_2\approx3.12$ with the minimal cost $J\approx45.9$. The trajectory calculations of Proposition~\ref{2-crowd} give us the same formulas \eqref{traj1} as in Cases~1.

Next we calculate the dual elements from Theorem~\ref{n-crowd} in this example. It follows from \eqref{e:a-q} that
\begin{eqnarray*}
\begin{cases}
-\frac{1}{\sqrt2}q^x_{11}(t)+\frac{1}{\sqrt2}q^x_{12}(t)=\dfrac{\lm\oa_1}{s_1}=\dfrac{3.12\lm}{6}=\dfrac{1.56}{3}=0.52,\\[2ex]
-\frac{1}{\sqrt2}q^x_{21}(t)+\frac{1}{\sqrt2}q^x_{22}(t)=\dfrac{\lm\oa_2}{s_2}=\dfrac{1.56\lm}{3}=0.52,
\end{cases}
\end{eqnarray*}
which gives us the following adjoint vector function when $\lm=1$:
\begin{eqnarray*}
q^x_{12}(t)-q^x_{11}(t)\approx0.74,\quad q^x_{22}(t)-q^x_{21}(t)\approx0.74\;\mbox{ for }\;t\in[0,6],
\end{eqnarray*}
while \eqref{e:qij} does not provide any extra information. We get an obvious solution of this system on $[0,6]$:
\begin{eqnarray*}
q^x_{11}(t)=0,\;q^x_{12}(t)=0.74,\;q^x_{21}(t)=0,\;q^x_{22}(t)=0.74.
\end{eqnarray*}
Remembering from the above that $\eta_{12}=\frac{9}{2}\oa_2=7.02$, conditions (5), (7), and (8) of Theorem~\ref{n-crowd} yield
\begin{eqnarray*}
p^x_1(t)=p^x_1(6)=(2.62,-2.62)-7.02\left(\frac{1}{\sqrt2},-\frac{1}{\sqrt2}\right)=(-2.34,2.34),\\
p^x_2(t)=p^x_2(6)=(-1.62,1.62)+7.02\left(\frac{1}{\sqrt2},-\frac{1}{\sqrt2}\right)=(3.34,-3.34).
\end{eqnarray*}
Finally, the measure $\gg$ can be calculated by (7) as
\begin{eqnarray*}
\gg([t,6])=q^x(t)-p^x(t)=(2.34,-1.6,-3.34, 4.08)\;\mbox{ for all }\;0\le t\le 6,
\end{eqnarray*}
which reflects the fact that the optimal sweeping motion hits the boundary of the state constraint \eqref{e:state constraints} at the initial time and stays there until the end of the process at $T=6$.}
\end{example}\vspace*{-0.1in}

It makes sense to discuss the situation in Example~\ref{ex1}, which provides hints for more general settings.\vspace*{-0.05in} 

\begin{remark}{\bf (analysis of the situation when participants are in contact).}\label{Rm1}
{\rm Although the optimal trajectories in both cases of Example~\ref{ex1} are the same, these cases reflect two largely different settings.\vspace*{-0.03in}

$\bullet$ In Case~1 the actual and spontaneous velocities have the same value, which is $(8.27,-8.27,8.27,-8.27)$. This shows that the normal cone $N_{C(t)}(\ox(t))$ is not active as $\eta_{12}=0$. It seems that to make the velocities of two participants coincide so that they can always maintain the constant distance $2R=6$ away from each other, we would think of the possibility of adjusting their initial different speeds such that $\oa_1s_1=\oa_2s_2$. And then the actual velocities will take the same value as the spontaneous velocities. This intuitive approach is very natural and seems to be reasonable at the first glance. However, the solution found in this case is proved {\em not to be optimal}. This makes a perfect sense since in reality the actual velocities and spontaneous velocities are unlikely to be the same when two participants are in contact.\vspace*{-0.03in}

$\bullet$ In Case~2 the actual velocities and the spontaneous velocities have different values, which are approximately $(8.27,-8.27,8.27,-8.27)$ and $(13.24,-13.24,3.31,-3.31)$, respectively. It is reasonable since the first participant, being farther to the exit than the second one is, tends to run faster based on his/her initial speed in the absence of the other. Thus this participant has to use more energy than the other uses with $\oa_1=2\oa_2\approx3.12$. The second participant does not need to be hurry since he/she is closer to the exit than the first one is. Getting in contact however, they both must adjust to the same actual velocity in order to maintain the constant distance from each other. In this case the normal cone plays a crucial role in the dynamical system with the active generating vector $\eta_{12}\nabla g(x)=7.02\left(-\frac{\sqrt2}{2},\frac{\sqrt2}{2}\right)$. This nonzero vector affects the motion of the two participants and causes them to use less energy comparing to the previous case. It is shown that the solution found in this case is actually an optimal solution.}
\end{remark}\vspace*{-0.1in}

The next two examples concern the motions when participants are {\em out of contact} at the initial time.\vspace*{-0.1in}

\begin{example}{\bf(participants are out of contact but have the same direction at the initial time).}\label{ex2}  {\rm Consider the optimal control problem in \eqref{e:state constraints}--\eqref{e:cost functional} with the following initial data:
\begin{eqnarray*}
n=2,\;T=6,\;s_1=6,\;s_2=3,\;x_{01}=(-60,60),\;x_{02}=(-48,48),\;R=3.
\end{eqnarray*}
In this case we have $\|\ox_2(0)-\ox_1(0)\|=12\sqrt2>6$, $\th_1(0)=\th_2(0)=135^\circ$, and $t_{12}>0$; see Figure~5.

\begin{center}
\begin{tikzpicture}
\draw [thick, red](-5,0)--(5,0);
\draw [thick, red] (0,0)--(0,6);
\draw [ultra thick, fill = orange] (-0.5,1) rectangle (0.5,0);
\draw node[below] at (0,0){Exit};
\draw node[below] at (0,-1){{\bf Figure 5}};
\draw [thick] (-4,4) circle (0.5);
\draw node[below] at (-4,4){$x_1$};
\draw [thick] (-2.5,2.5) circle (0.5);
\draw node[below] at (-2.5,2.5){$x_2$};
\draw [ultra thick, blue] [->] (-4,4) -- (-3.2,3.2);
\draw [ultra thick, blue] [->] (-2.5,2.5) -- (-2,2);
\coordinate[label=above left:] (A) at (-4,4);
\coordinate[label=above left:] (B) at (-5,5);
\coordinate[label=above left:] (C) at (-2.5,4);
\draw[->,purple] (A)--(B);
\draw[->,purple] (A)--(C);
\draw[dashed] (-4,4)--(-2.5,2.5);
\markangle{B}{A}{C}{$135^\circ$}{2}
\end{tikzpicture}
\end{center}
It follows from \eqref{e:t12} that $\eta_{12}>0$ and hence $\oa_1=2\oa_2$. Furthermore, we get from \eqref{e:eta cal} and \eqref{e:dir21} that $\cos\th_{21}=\frac{\sqrt2}{2}$ and $\sin\th_{21}=-\frac{\sqrt2}{2}$. Thus $\eta_{12}=\frac{1}{2}(6\oa_1-3\oa_2)=\frac{9}{2}\oa_2$, and the optimal trajectories after the contact time are calculated by
\begin{eqnarray*}
\left\{
\begin{array}{ll}
\ox_1(t)=\left(-60+\frac{15\sqrt2}{4}\oa_2t+\frac{9\sqrt2}{4}\oa_2t_{12},60-\frac{15\sqrt2}{4}\oa_2t-\frac{9\sqrt2}{4}\oa_2t_{12}\right),\\[2ex]
\ox_2(t)=\left(-48+\frac{15\sqrt2}{4}\oa_2t-\frac{9\sqrt2}{4}\oa_2t_{12},48-\frac{15\sqrt2}{4}\oa_2t+\frac{9\sqrt2}{4}\oa_2t_{12}\right)
\end{array}\right.
\;\mbox{ for }\;t\in[t_{12},6].
\end{eqnarray*}
Using \eqref{e:t12} gives us $\eta_{12}t_{12}=6\sqrt2-3$, and so $\oa_2t_{12}=\dfrac{4\sqrt2}{3}-\dfrac{2}{3}$. As a consequence, we get
\begin{eqnarray*}
\left\{\begin{array}{ll}
\ox_1(t)=\left(-54-\frac{3\sqrt2}{2}+\frac{15\sqrt2}{4}\oa_2t,54+\frac{3\sqrt2}{2}-\frac{15\sqrt2}{4}\oa_2t\right),\\[2ex]
\ox_2(t)=\left(-54+\frac{3\sqrt2}{2}+\frac{15\sqrt2}{4}\oa_2t,54-\frac{3\sqrt2}{2}-\frac{15\sqrt2}{4}\oa_2t\right)
\end{array}\right.
\;\mbox{ for }\;t\in[t_{12},6].
\end{eqnarray*}
Substituting the values of $\ox_1(6)$ and $\ox_2(6)$ into the cost functional \eqref{e:cost functional} shows that
\begin{eqnarray*}
J[x,a]=2040\oa_2^2-2\cdot2430\sqrt2\oa_2+\left(54+\frac{3\sqrt2}{2}\right)^2+\left(-54+\frac{3\sqrt2}{2}\right)^2.
\end{eqnarray*}
The minimum point here is $\oa_2=\frac{2430\sqrt2}{2040}\approx 1.68$, which allows us to compute the optimal solution as follows: $(\oa_1,\oa_2)=(3.36,1.68)$,
\begin{eqnarray*}
\ox_1(t)=\left\{\begin{array}{ll}
(14.29t-60,-14.29t+60)\;\mbox{ for }\;t\in[0,0.72),\\[1ex]
\left(8.91t-54-\frac{3\sqrt2}{2},-8.91t+54+\frac{3\sqrt2}{2}\right)\;\mbox{ for }\;t\in[0.72,6];
\end{array}\right.
\end{eqnarray*}
\begin{eqnarray*}
\ox_2(t)=\left\{\begin{array}{ll}
(3.56t-48,-3.56t+48)\;\mbox{ for }\;t\in[0,0.72),\\[1ex]
\left(8.93t-54+\frac{3\sqrt2}{2},-8.93t+54-\frac{3\sqrt2}{2}\right)\;\mbox{ for }\;t\in[0.72,6].
\end{array}\right.
\end{eqnarray*}
Thus the participants adjust to the same velocity and keep it till the end of the process; see Figure~6.

\begin{center}
\begin{tikzpicture}
\draw [thick, red](-5,0)--(5,0);
\draw [thick, red] (0,0)--(0,6);
\draw [ultra thick, fill = orange] (-0.5,1) rectangle (0.5,0);
\draw node[below] at (0,0){Exit};
\draw node[below] at (0,-1){{\bf Figure 6}};
\draw [thick] (-3,3) circle (0.5);
\draw node[below] at (-3,3){$x_1$};
\draw [thick] (-2.3,2.3) circle (0.5);
\draw node[below] at (-2.3,2.3){$x_2$};
\draw [ultra thick, blue] [->] (-3,3) -- (-2.5,2.5);
\draw [ultra thick, blue] [->] (-2.3,2.3) -- (-1.8,1.8);
\draw[dashed] (-3,3)--(-2.3,2.3);
\coordinate[label=above left:] (A) at (-3,3);
\coordinate[label=above left:] (B) at (-4,4);
\coordinate[label=above left:] (C) at (-2,3);
\draw[->,purple] (A)--(B);
\draw[->,purple] (A)--(C);
\markangle{B}{A}{C}{$135^\circ$}{2}
\end{tikzpicture}
\end{center}
We next compute the dual elements in this model. It follows from \eqref{e:a-q} that
\begin{eqnarray*}
\begin{cases}
-\frac{1}{\sqrt2}q^x_{11}(t)+\frac{1}{\sqrt2}q^x_{12}(t)=\dfrac{\lm\oa_1}{s_1}=\dfrac{3.36\lm}{6}=0.56,\\[2ex]
-\frac{1}{\sqrt2}q^x_{21}(t)+\frac{1}{\sqrt2}q^x_{22}(t)=\dfrac{\lm\oa_2}{s_2}=\dfrac{1.68\lm}{3}=0.56,
\end{cases}
\end{eqnarray*}
which leads us to the relationships
\begin{eqnarray*}
q^x_{12}(t)-q^x_{11}(t)\approx0.79,\;q^x_{22}(t)-q^x_{21}(t)\approx0.79\;\mbox{ on }\;[0,6]
\end{eqnarray*}
provided that $\lm=1$. Similarly to Example~\ref{ex1}, we choose one obvious solution to these equations
\begin{eqnarray*}
q^x_{11}(t)=0,\;q^x_{12}(t)=0.79\,\;q^x_{21}(t)=0,\;q^x_{22}(t)=0.79\;\mbox{ on }\;[0,6]
\end{eqnarray*}
and calculate $\eta_{12}=\frac{9}{2}\oa_2=7.56$ by \eqref{e:t12}. Further, we get from (5), (6), and (8) of Proposition~\ref{2-crowd} that
\begin{eqnarray*}
\begin{cases}
p^x_1(t)=p^x_1(6)=(2.66,-2.66)-7.56\left(\frac{1}{\sqrt2},-\frac{1}{\sqrt2}\right)=(-2.69,2.69),\\[2ex]
p^x_2(t)=p^x_2(6)=(-1.7,1.7)+7.56\left(\frac{1}{\sqrt2},-\frac{1}{\sqrt2}\right)=(3.65,-3.65)
\end{cases}
\end{eqnarray*}
on $[0,6]$, which yields by (7) the measure expression
\begin{eqnarray*}
\gg([t,6])=q^x(t)-p^x(t)=(2.69,-1.9,-3.65,4.44)\;\mbox{ for }\;0.72\le t\le 6.
\end{eqnarray*}
The latter reflects the fact that the optimal sweeping motion hits the boundary of the state constraints at time $t_{12}=0.72$ and stays there until the end of the process at $T=6$.}
\end{example}\vspace*{-0.1in}

The last example illustrates the situation when two participants are out of contact while having different directions at the beginning of the dynamic process.\vspace*{-0.1in}

\begin{example}{\bf(participants are out of contact and have different directions at the initial time).}\label{ex3}
{\rm Consider the optimal control problem  in \eqref{e:state constraints}--\eqref{e:cost functional} with the following initial data:
\begin{eqnarray*}
n=2,\;T=6,\;s_1=6,\;s_2=3,\;x_{01}=\left(-60,60\right),\;x_{02}=(-48,54),\;R=3.
\end{eqnarray*}
In this case we have $\|\ox_2(0)-\ox_1(0)\|=6\sqrt5>6, \th_1(0)=135^\circ$, and $\th_2(0)=131.63^\circ$ as shown in Figure~7. Thus $t_{12}>0$ and $\eta_{12}>0$ due to \eqref{e:t12}. It follows from \eqref{e:dir21} that $\cos\th_{12}=\dfrac{2}{\sqrt5}$ and $\sin\th_{12}=-\dfrac{1}{\sqrt5}$. On the other hand, we deduce from \eqref{e:eta cal} and \eqref{e:eta cal1} that
\begin{eqnarray}\label{e:eta cal2}
(\cos\th_{12})\eta_{12}=-\frac{12}{\sqrt{145}}\oa_2+\frac{3\sqrt2}{2}\oa_1,\quad
(\sin\th_{12})\eta_{12}=\frac{27}{2\sqrt{145}}\oa_2-\frac{3\sqrt2}{2}\oa_1,
\end{eqnarray}
which implies that $\oa_2=\frac{\sqrt{290}}{10}\oa_1$ and $\eta_{12}=\frac{3\sqrt{10}}{20}\oa_1$.

\begin{center}
\begin{tikzpicture}
\draw [thick, red](-5,0)--(5,0);
\draw [thick, red] (0,0)--(0,6);
\draw [ultra thick, fill = orange] (-0.5,1) rectangle (0.5,0);
\draw node[below] at (0,0){Exit};
\draw node[below] at (0,-1){{\bf Figure 7}};
\draw [thick] (-4,4) circle (0.5);
\draw node[below] at (-4,4){$x_1$};
\draw [thick] (-3,3.375) circle (0.5);
\draw node[below] at (-3,3.375){$x_2$};
\draw [ultra thick, blue] [->] (-4,4) -- (-3.2,3.2);
\draw [ultra thick, blue] [->] (-3,3.375) -- (-2.5,2.8125);
\coordinate[label=above left:] (A) at (-4,4);
\coordinate[label=above left:] (B) at (-5,5);
\coordinate[label=above left:] (C) at (-3,4);
\draw[->,purple] (A)--(B);
\draw[->,purple] (A)--(C);
\markangle{B}{A}{C}{$135^\circ$}{2}
\coordinate[label=above left:] (A') at (-3,3.375);
\coordinate[label=above left:] (B') at (-3.5,3.94);
\coordinate[label=above left:] (C') at (-2,3.375);
\draw[->,purple] (A')--(B');
\draw[->,purple] (A')--(C');
\markangle{B'}{A'}{C'}{$131.63^\circ$}{2}
\end{tikzpicture}
\end{center}
The trajectories before and after the contact time are calculated as follows:
\begin{eqnarray*}
\begin{cases}
\ox_1(t)=(-60+3\sqrt2\oa_1t,60-3\sqrt2\oa_1t),\\[1ex]
\ox_2(t)=\left(-48+\frac{12\sqrt2}{5}\oa_1t, 54-\frac{27\sqrt2}{10}\oa_1t\right)
\end{cases}
\;\mbox{ for }\;t\in[0,t_{12});
\end{eqnarray*}
\begin{eqnarray*}
\begin{cases}
\ox_1(t)=\left(-54-\frac{6}{\sqrt5}+\frac{27\sqrt2}{10}\oa_1t, 57+\frac{3}{\sqrt5}-\frac{57\sqrt2}{20}\oa_1t\right),\\[1ex]
\ox_2(t)=\left(-54+\frac{6}{\sqrt5}+\frac{27\sqrt2}{10}\oa_1t, 57-\frac{3}{\sqrt5}-\frac{57\sqrt2}{20}\oa_1t\right)
\end{cases}
\;\mbox{ for }\;t\in[t_{12},6].
\end{eqnarray*}
Substituting the values of $\ox_1(6)$ and $\ox_2(6)$ into the cost functional \eqref{e:cost functional} gives us the function of $\oa_1$: 
\begin{eqnarray*}
\begin{aligned}
J[\ox,\oa]&=\frac{1}{2}\bigg[\left(-54-\frac{6}{\sqrt5}+\frac{81\sqrt2}{5}\oa_1\right)^2+\left(57+\frac{3}{\sqrt5}-\frac{171\sqrt2}{10}\oa_1\right)^2+\left(-54+\frac{6}{\sqrt5}+\frac{81\sqrt2}{5}\oa_1\right)^2\\
&+\left(57-\frac{3}{\sqrt5}-\frac{171\sqrt2}{10}\oa_1\right)^2\bigg]+11.7\oa_1^2\\
&=1121.4\oa_1^2-3699\oa_1+(54+6/\sqrt5)^2+(-54+6/\sqrt5)^2+(57+3/\sqrt5)^2+(57-3/\sqrt5)^2,
\end{aligned}
\end{eqnarray*}
which clearly attains its minimum at $\oa_1=\frac{3699}{2\times1121.4}\approx1.65$ with $\oa_2=\frac{290}{10}\oa_1\approx 2.8$. Then the formulas of Proposition~\ref{2-crowd} for optimal trajectories give us the expressions
\begin{eqnarray*}
\ox_1(t)=(-60+7t,60-7t),\;\ox_2(t)=\left(-48+5.6t,54-6.3t\right)\;\mbox{ on }\;[0,4.74);
\end{eqnarray*}
\begin{eqnarray*}
\begin{cases}
\ox_1(t)=\left(-54-\frac{6}{\sqrt5}+6.3t, 57+\frac{3}{\sqrt5}-6.65t\right),\\[1ex]
\ox_2(t)=\left(-54+\frac{6}{\sqrt5}+6.3t, 57-\frac{3}{\sqrt5}-6.65t\right)
\end{cases}
\;\mbox{ for }\;t\in[4.74,6].
\end{eqnarray*}
As illustrated by Figure~8, the two participants switch to the new directions and the new velocities after being in contact and then maintain their new velocities till the end of the process.

\begin{center}
\begin{tikzpicture}
\draw [thick, red](-5,0)--(5,0);
\draw [thick, red] (0,0)--(0,6);
\draw [ultra thick, fill = orange] (-0.5,1) rectangle (0.5,0);
draw node[below] at (0,0){Exit};
\draw node[below] at (0,-1){{\bf Figure 8}};
\draw [thick] (-1.4,1.5) circle (0.5);
\draw node[below] at (-2.1,2.2){$x_1$};
\draw [thick] (-2.1,2.2) circle (0.5);
\draw node[below] at (-1.4,1.5){$x_2$};
\draw [ultra thick, blue] [->] (-1.4,1.5) -- (-0.9,0.7);
\draw [ultra thick, blue] [->] (-2.1,2.2) -- (-1.5,1.24);
\coordinate[label=above left:] (A) at (-2.1,2.2);
\coordinate[label=above left:] (B) at (-3,3.16);
\coordinate[label=above left:] (C) at (-1.1,2.2);
\end{tikzpicture}
\end{center}
Next we compute the dual elements in this example. It follows from \eqref{e:a-q} that
\begin{eqnarray*}
\left\{\begin{array}{ll}
-\frac{1}{\sqrt2}q^x_{11}(t)+\frac{1}{\sqrt2}q^x_{12}(t)=\dfrac{\lm\oa_1}{s_1}=\dfrac{1.65}{6}=0.275,\\[2ex]
-\frac{8}{\sqrt{145}}q^x_{21}(t)+\frac{9}{\sqrt{145}}q^x_{22}(t)=\dfrac{\lm\oa_2}{s_2}=\dfrac{2.8}{3}=0.93
\end{array}\right.
\end{eqnarray*}
on $[0,6]$ provided that $\lm=1$. Since $\eta_{12}(t)=\frac{3\sqrt{10}}{20}\oa_1\approx 0.78$ for all $t\in[4.74 6]$, we deduce from \eqref{e:qij} that
\begin{eqnarray*}
\frac{12}{\sqrt5}\big(q^x_{21}(t)-q^x_{11}(t)\big)-\frac{6}{\sqrt5}\big(q^x_{22}(t)-q^x_{12}(t)\big)=0.
\end{eqnarray*}
Combining it with the equations above gives us the linear system
\begin{eqnarray*}
\begin{cases}
-q^x_{11}(t)+q^x_{12}(t)=0.39,\\
-8q^x_{21}(t)+9q^x_{22}(t)=11.2,\\
q^x_{22}(t)-q^x_{12}(t)=2\big(q^x_{21}(t)-q^x_{11}(t)\big)
\end{cases}\;
\mbox{ for all }\;t\in[4.74,6],
\end{eqnarray*}
which has infinitely many solutions. Similarly to the previous examples, we pick any particular one and then find the other dual elements $p^x_1(t),p^x_2(t),$ and $\gg([t,6])$ for $t\in[4.74,6]$ from the conditions of Theorem~\ref{n-crowd}. We would also see that the measure $\gg([t,6])$ reduces to a positive constant on $[4.74,6]$, which reflects the fact that the optimal motion hits the boundary of the state constraint and stays there till the end of the process at $t=6$. It is worth mentioning that in this case we do not get the equality $s_2\oa_1=s_1\oa_2$, since the two participants do not have the same direction at the contact time. Thus equation \eqref{e:a-q} does not provide anymore useful information about the link between $\oa_1$ and $\oa_2$. However, such a relationship can be found by using \eqref{e:eta cal2} that is valid for the two participants having different directions at the initial time. For brevity, we skip the further standard calculations.}
\end{example}\vspace*{-0.2in}

\section{Concluding Remarks}\vspace*{-0.1in}

This paper demonstrates that the necessary optimality conditions obtained for a rather general nonconvex version of the controlled sweeping process with prox-regular moving sets are instrumental to control and optimize a practical planar crowd motion model formalized as a sweeping process of the aforementioned type. Deriving necessary optimality conditions for optimal motions with finitely many participants, we present their complete analytical realization in the case of two participants. The major question remains on developing efficient numerical algorithms to solve the obtained systems of optimality conditions in crowd motion models with many participants. This is a challenging issue for our future research. Furthermore, it seems possible to apply the developed variational machinery to optimization problems in robotics, hysteresis, and systems of engineering design that are modeled as controlled sweeping processes. These topics are also among our subsequent research goals. 
\end{document}